%% file: agt-5-46.tex
\documentclass{gtart_h}
\input agtout

\lognumber{46}
\volumenumber{5}
\volumeyear{2005}
\papernumber{46}
\pagenumbers{1141}{1172}
\received{25 November 2004} 
\revised{27 August 2005}
\accepted{29 August 2005}
\published{18 September 2005}

\usepackage[all]{xy}
\SelectTips{cm}{}
\usepackage{amssymb}

\newtheorem{thm}{Theorem}[section]
\newtheorem{cor}[thm]{Corollary}
\newtheorem{lem}[thm]{Lemma}
\newtheorem{prop}[thm]{Proposition}
\theoremstyle{definition}
\newtheorem{defn}{Definition}
\newtheorem{rem}[thm]{Remark}
\newtheorem{exmp}[thm]{Example}
\newcommand{\s}{{\mathcal S}}

\begin{document}

\shorttitle{Nullification functors and the homotopy type of classifying spaces}
\title{Nullification functors and the homotopy type\\of the classifying space
for proper bundles}

\author{Ram\'on J. Flores }
\coverauthors{Ram\noexpand\'on J. Flores }
\asciiauthors{Ramon J. Flores }

\address{Departamento de Matem\'aticas, Universidad Aut\'onoma 
de Barcelona\\E--08193 Bellaterra, Spain} 
\asciiaddress{Departamento de Matematicas, Universidad Autonoma 
de Barcelona\\E-08193 Bellaterra, Spain}

\email{ramonj@mat.uab.es}

\urladdr{http://mat.uab.es/~ramonj}
\asciiurl{http://mat.uab.es/ ramonj}

\begin{abstract}
Let $G$ be a discrete group for which the classifying space for
proper $G$-actions is finite-dimensional. We find a space $W$ such
that for any such $G$, the classifying space
$\underline{\textrm{B}}G$ for proper $G$-bundles has the homotopy
type of the $W$-nullification of $\textrm{B}G$. We use this to
deduce some results concerning $\underline{\textrm{B}}G$ and in
some cases where there is a good model for
$\underline{\textrm{B}}G$ we obtain information about the
$\textrm{B}\mathbb{Z}/p$-nullification of $\textrm{B}G$.
\end{abstract}

\asciiabstract{Let G be a discrete group for which the classifying space for
proper G-actions is finite-dimensional.  We find a space W such that
for any such G, the classifying space PBG for proper G-bundles has the
homotopy type of the W-nullification of BG. We use this to deduce some
results concerning PBG and in some cases where there is a good model
for PBG we obtain information about the BZ/p-nullification of BG.}


\maketitle

\section{Introduction}

Let $G$ be a discrete group. We will say that a
\emph{G}-CW-complex $X$ is \emph{proper} if the isotropy groups of
the action are finite.

In 1971, J.P. Serre introduced in \cite{Serre71} the ``classifying
space for proper actions", that can be described as the unique
proper \emph{G}-CW-complex $\underline{\textrm{E}}G$, up to
$G$-homotopy, that enjoys the following universal property:

``If $X$ is another proper \emph{G}-CW-complex, there exists a
$G$-map $X\longrightarrow\underline{\textrm{E}}G$ which is unique
up to $G$-homotopy".

The space $\underline{\textrm{E}}G$ appears as the principal new
feature in the reformulation of the Baum-Connes conjecture stated
in \cite{Baum94} by Baum-Connes-Higson. The conjecture, partially
solved, asks if for a locally compact, Hausdorff and second
countable group $G$, the assembly map from the Kasparov
$K$-homology groups $K_j^G(\underline{\textrm{E}}G)$ to the
$C^*$-algebra $K$-theory groups $K_j(C_r^*(G))$ is an isomorphism
for $j=0,1$. The great amount of research that has emerged around
this subject has led to a growing interest in the theory of proper
actions.

An important part of the efforts carried out in this direction has
been devoted to understand the relationship between the algebraic
structure of $G$ and the homotopy-theoretic properties of
$\underline{\textrm{E}}G$ and its quotient space
$\underline{\textrm{E}}G/G$, which is currently denoted by
$\underline{\textrm{B}}G$. Probably the greatest success has been
reached by interpreting correctly finiteness group-theoretic
conditions over $G$ in order to build models of
$\underline{\textrm{E}}G$ enjoying various types of finiteness
conditions. See (\cite{Lueck04}, section 5) for an excellent
survey on this topic.

In the same way that happens with classical $G$-actions, the
importance of $\underline{\textrm{E}}G$ and in particular of
$\underline{\textrm{B}}G$ does not come only because they reflect
geometrically the algebraic properties of the group $G$, but
because the importance of these spaces in the theory of
$G$-bundles. Baum-Connes-Higson already pointed out that
$\underline{\textrm{B}}G$ classifies proper $G$-bundles (see
\cite{Baum94}, definition 8), and they described how to obtain
them by making pullback on maps
$X\longrightarrow\underline{\textrm{B}}G$, a method of a clear
classical flavor. Moreover, if $G$ is a group for which
Baum-Connes conjecture holds, the knowledge of the homotopical
structure of $\textrm{\underline{B}}G$ plays a r\^ole in the
computation of the $K$-theory of $C^*(G)$, and in general gives
information about the integral homology of $G$. See
(\cite{Lueck04}, section 8) for details about these results.

The most important attempt made so far to understand the homotopy
type and properties of $\underline{\textrm{B}}G$ is the paper of
Leary-Nucinkis \cite{Leary01}. In it, the authors prove that for
every CW-complex $X$ there exists a discrete group $G_X$ such that
$\underline{\textrm{B}}G_X$ is homeomorphic to $X$. This
``Kan-Thurston-like" result is proved using essentially tools of
the theory of graph of groups. As a by-product, they obtain a
precise description of the fundamental group of
$\underline{\textrm{B}}G$ and a construction of
$\underline{\textrm{B}}G$ for some subgroups of right-angled
Coxeter groups.

Although these results have been very useful for us (particularly
the formula for the fundamental group), our approach to the
homotopy type of $\underline{\textrm{B}}G$ has been different, and
has been carried out with pure homotopy-theoretic tools. Our idea
is to find a functor $F$ in the topological category that
transforms models of $\textrm{B}G$ on spaces that are homotopy
equivalent to models of $\underline{\textrm{B}}G$. This functor
have enough good properties in order to read information about
$\underline{\textrm{B}}G$ from $\textrm{B}G$ and vice versa.

The appropriate functor turns to be a nullification; a tool that was
introduced by Bousfield in \cite{Bousfield94} in order to study
periodic phenomena in unstable homotopy (in fact, he called it
``periodization"), and that has been widely used since then.  The
utility of this functor in this context comes from the fact that it
will allow us to apply all the machinery of localization developed in
the 90's by Bousfield himself \cite{Bousfield94, Bousfield97},
Dror-Farjoun \cite{Dror-Farjoun95}, Chach\'olski
\cite{Chacholski96} and others.

Our main result is the following:

\medskip

{\bf Theorem \ref{bbarnull}}\qua\textsl{Let $G$ be a
discrete group such that there exists a finite-dimensional model
for $\underline{\emph{B}}G$. Let $W_{\infty}$ denote
$\bigvee\emph{B}\mathbb{Z}/p$, where the wedge is indexed by all
primes. Then we have a homotopy equivalence
$\mathbf{P}_{W_{\infty}}\emph{B}G\simeq\underline{\emph{B}}G$,
where $\mathbf{P}_{W_{\infty}}$ denotes the
$W_{\infty}$-nullification functor.}

\medskip

Observe that the condition concerning the existence of a model of
$\underline{\textrm{B}}G$ with finiteness conditions is not too
restrictive, in account of the great quantity of groups that have
recently appeared for which these conditions hold.

Now we describe in more detail the contents of each of the
sections of the paper, and in particular we will comment a little
bit the consequences of the main theorem.

In section 2 we review the needed background and results
concerning proper actions, and we construct a particular model for
the classifying space for families $E_{\mathcal{F}}G$ that will be
very useful in the rest of the work; as a fundamental by-product,
we also obtain an appropriate model for
$\textrm{B}_{\mathcal{F}}G$.

Section 3 constitutes the bulk of the work, because it is devoted
to the proof of the main theorem we stated before. The technique
is the following: we apply the functor $\mathbf{P}_{W_{\infty}}$
to a suitable model of $\textrm{B}G$, and we obtain that it is
homotopy equivalent to the $W_{\infty}$-nullification of the nerve
of some small category that only depends on $G$. This nerve turns
out to be the model built in the previous section for
$\underline{\textrm{B}}G$, and we finish by checking that in the
conditions of the theorem $\underline{\textrm{B}}G$ is
$W_{\infty}$-null.

The rest of the paper is devoted to take out some consequences of
the main theorem. So, in section 4 we describe the behaviour of
the functor $\textrm{\underline{B}}$ with respect to various
fundamental constructions in homotopy theory, namely products,
wedges or colimits. Moreover, we identify in some cases the
universal cover of $\textrm{\underline{B}}G$ and we obtain some
conditions about preservation of fibrations under passing to
classifying spaces.

We begin the following section with a short new proof of the
well-known fact that if $G$ is a locally finite group which
cardinal is smaller than $\aleph_{\omega}$,
$\underline{\textrm{B}}G$ is contractible (indeed, the statement
is true for all locally finite groups, see \cite{DKLT02}). Later,
we treat the case of groups for which the normalizer condition
holds, an ample class of discrete groups that include, for
example, all the nilpotent groups. We prove that if a group $G$ in
this class admits a finite-dimensional model for
$\underline{\textrm{B}}G$, then
$\underline{\textrm{B}}G\simeq\textrm{B}H$ for some quotient group
$H$ that we identify. In particular, in this case
$\underline{\textrm{B}}G$ is nilpotent as a space if $G$ is
nilpotent as a group. We finish this paragraph by studying the
$\textrm{B}\mathbb{Z}/p$-nullification of classifying spaces of
supersoluble groups.

In section 6 we take the opposite point of view, showing that the
main theorem can give information in the two directions. More
concretely, we focus our attention on groups of isometries of the
real plane, and taking profit of well-known geometric properties
of them we obtain via $\textrm{\underline{B}}G$ a lot of
information on the $\textrm{B}\mathbb{Z}/p$-nullification of their
classifying spaces.

Last section is devoted to the study of a canonical map
$\textrm{B}G\longrightarrow\textrm{B}_{\mathcal{F}}G$ (defined
previously) that always relates the classical and proper
classifying spaces. More concretely, we prove that the homotopy
fiber of that map can be described as a homotopy colimit of
classifying spaces of groups of $\mathcal{F}$ over a contractible
category. We finish by proving a technical and interesting
statement that appears in the proof and concerns the localization
of a comma category.

\textbf{Acknowledgements}\qua I wish to thank Carles Broto, for turning
my attention to proper actions and the world of Geometric Group
Theory, and for all the time we have spent discussing about these
topics. The results about crystallographic groups were motivated
by a suggestion of Ian Leary, whom I acknowledge his interest in
my work. I am also grateful to Emmanuel Dror-Farjoun, who pointed
out some very useful observations that enriched the results of
last section. Finally, I would like to thank the Institute
Galil\'ee, Universit\'e Paris XIII, for their hospitality in the
seven months I spent there, and the referee, for the detailed
report and numerous suggestions.
The author is partially supported by MCYT grant BFM2001-2035.

\section{Useful models for $\textrm{E}_{\mathcal{F}}G$ and $\textrm{B}_{\mathcal{F}}G$}

We will begin by recalling the definition of the classifying space
for families. Further information about proper actions and their
classifying spaces can be found in \cite{Baum94}, \cite{Leary01},
\cite{MislinV01}, or \cite{Dieck87}.

\begin{defn}

Suppose $\mathcal{F}$ is a family of subgroups of a discrete group
$G$ that is closed under conjugation and taking subgroups. We will
say that a $G$-CW-complex $Y$ is a \emph{model} for
$\emph{E}_{\mathcal{F}}G$ if the isotropy group of each point
belongs to $\mathcal{F}$ and for each $H\in\mathcal{F}$, the
fixed-point space $Y^H$ is contractible.

\end{defn}

The $G$-space $\textrm{E}_{\mathcal{F}}G$ is characterized by the
following universal property:

\begin{prop}

If $X$ is a model for $\emph{E}_{\mathcal{F}}G$, then for each
$G$-CW-complex $Y$ whose isotropy groups lie in $\mathcal{F}$
there is a map $Y\longrightarrow X$ which is unique up to
$G$-homotopy. Moreover, two models for $\emph{E}_{\mathcal{F}}G$
are always $G$-homotopy equivalent. Conversely, if $X$ is a
$G$-CW-complex for which this universal property holds, $X$ is a
model for $\emph{E}_{\mathcal{F}}G$.

\end{prop}

The quotient space $\textrm{E}_{\mathcal{F}}G/G$ is usually
denoted by $\textrm{B}_{\mathcal{F}}G$. As
$\textrm{E}_{\mathcal{F}}G$ is unique up to $G$-homotopy
equivalence, $\textrm{B}_{\mathcal{F}}G$ is unique up to homotopy
equivalence; see (\cite{Leary01}, section 2) for a review of the
main properties of $\textrm{B}_{\mathcal{F}}G$.

If there is no explicit mention against it, we will suppose from
here that $\mathcal{F}$ is the family of finite subgroups of $G$.
In this case, it is standard to denote $\textrm{E}_{\mathcal{F}}G$
and $\textrm{B}_{\mathcal{F}}G$ by $\underline{\textrm{E}}G$ and
$\underline{\textrm{B}}G$, respectively. In particular, it can be
seen that every group homomorphism $G\longrightarrow G'$ induces a
$G$-map
$\underline{\textrm{E}}G\longrightarrow\underline{\textrm{E}}G'$
(respectively a map
$\underline{\textrm{B}}G\longrightarrow\underline{\textrm{B}}G'$)
which is unique up to $G$-homotopy (respectively up to homotopy).

A number of geometric constructions for
$\textrm{E}_{\mathcal{F}}G$ and its quotient space is available in
the literature (see section 4 of \cite{Lueck00} for a description
of some of them), but to describe the relationship between
$\textrm{B}G$ and $\textrm{B}_{\mathcal{F}}G$ we will need to
build these spaces as nerves of small categories. This is the main
goal of this section.

Recall first the definition of the \emph{orbit category}
associated to a group $G$ and a family of subgroups $\mathcal{F}$:
the objects are the homogeneous spaces $G/H$, with
$H\in\mathcal{F}$, and the morphisms are the $G$-maps. It is not
hard to see that there is a bijective map
$$\textrm{Mor}(G/K ,G/H)=\{g\in G\textrm{ }|\textrm{ }g^{-1}Kg\subseteq H\}/H$$
given by $f\longrightarrow f(eK)$, where $e$ is the identity
element of $G$. The key definition we need for building the
desired model of $\textrm{E}_{\mathcal{F}}G$ is the following (see
\cite{Dwyer97}, section 2, for details):
\begin{defn}
Let $\mathcal{D}$ be a small category, $\mathbf{Cat}$ the category
where the objects are the small categories and whose morphisms are
functors, and
 $f:\mathcal{D}\rightarrow\mathbf{Cat}$ a functor. The \emph{Grothendieck
 construction} $\mathbf{Gr}(f)$ associated to $f$
 is defined as the category whose objects are the pairs $(d,x)$,
 with $x\in\mathcal{D}$ and $x\in f(d)$, and where a morphism $(d,x)\rightarrow
 (d',x')$ is a pair $(u,v)$ where $u:d\rightarrow d'$ is a
 morphism in $\mathcal{D}$ and $v:f(u)(x)\rightarrow x'$ is a
 morphism in $f(d')$. The composition is made in the obvious way.
 \end{defn}

The main feature of this construction, due to Thomason, is the
following:

\begin{thm}
Let $\mathcal{D}$ be a small category,
$F:\mathcal{D}\longrightarrow\mathbf{Cat}$ a functor, and
$\mathbf{Gr}(F)$ the Grothendieck construction of $F$. Then there
exists a natural weak homotopy equivalence:
$$\emph{N}(\mathbf{Gr}(F))\simeq\emph{hocolim N}(F)$$ where
$\emph{N}$ denotes the nerve.

\label{Thomason}
\end{thm}

\begin{proof}

See \cite{Thomason79}, 1.2.
\end{proof}

Now we can describe our model of the universal space
$\textrm{E}_{\mathcal{F}}G$.

\begin{prop}
Let $G$ be a discrete group. Consider the functor
$$R:\mathbf{O}_{\mathcal{F}}\longrightarrow\mathbf{Cat}$$
that sends every homogeneous space $G/H$ to the category $G/H$
(whose objects are the elements of $G/H$ and there is only
identity morphisms), and the morphisms to the obvious functors. In
these conditions, we have that the nerve of the Grothendieck
construction of $f$ is a model for $\emph{E}_{\mathcal{F}}G$.
\label{modelo}
\end{prop}

\begin{proof}

For convenience, we will denote $X=|\textrm{N}(\mathbf{Gr}(R))|$.
This space  has a natural action of $G$ given by the left action
of $G$ in every homogeneous space $G/H$. We will prove firstly
that for every $x\in X$ the isotropy group $G_x$ of $x$ belongs to
$\mathcal{F}$.

The action of $G$ over $X$ is simplicial and is induced from the
action of $G$ over the homogeneous spaces, so by definition of
nerve we will only need to study the action of $G$ over the
vertices of $X$. Hence, let $(G/H,a)$ be the pair associated with
the vertex $x$. It is clear that
$$G_x=\{g\in G\textrm{ }|\textrm{ }gaH=aH\}=\{g\in G\textrm{
}|\textrm{ }\exists h\in H\textrm{ s. t. } gah=a\}$$ and this
amounts to say that $g\in aHa^{-1}$. So, $G_x=aHa^{-1}$, that
belongs to $\mathcal{F}$ because $H$ does.

Now we will see that for every $K\in\mathcal{F}$, the set of fixed
points $X^K$ is contractible. Now, a point on a simplex is fixed
by a subgroup if and only if all vertices are fixed, because the
action of $G$ is defined over the vertices and then extended to
the rest of the simplices by linearity.

So, consider a vertex $x\in X$. This point is fixed by the action
of $K$ if, given the pair $(G/H,aH)$ associated to it, we have
that, for every $k\in K$, $kaH=aH$. Thus, we see that
$$X^K=\bigcup_{H\in\mathcal{F}}\{(G/H,aH)\textrm{ }|\textrm{ }a^{-1}Ka\subseteq H\}/H.$$
So, for every element $(G/H, aH) \in X^K$ there exists one and
only one morphism $(G/K,eK)\longrightarrow (G/H,aH)$, or in other
words, $X^K$ can be identified with the nerve of the undercategory
(in the sense of \cite{MacLane71}, II.6) associated to the element
$(G/K,eK)$ of the category $\mathbf{Gr}(R)$. This nerve is
contractible, and then $X^K$ is contractible too.  So we are done.
 \end{proof}

Note that if we consider the action of $G$ over $\textbf{Gr}(R)$
via functors, the objects of $\textbf{Gr}(R)/G$ are the
homogeneous spaces $G/H$, $H<G$ finite, and the morphisms are the
$G$-maps. So, this quotient category is identified in a natural
way with the orbit category $\mathbf{O}_{\mathcal{F}}$, and in
particular $\textrm{N}(\mathbf{O}_{\mathcal{F}})$ is a model for
$\textrm{B}_{\mathcal{F}}G$.

\begin{rem}

Observe that if we take the realizations of the nerves of
$\textbf{\emph{Gr}}(R)$ and $\mathbf{O}_{\mathcal{F}}$, we obtain
models of $\emph{E}_{\mathcal{F}}G$ and $\emph{B}_{\mathcal{F}}G$
in the category of (topological) spaces.

\end{rem}

\begin{rem}
The idea of the construction comes from (\cite{Arone01}, section
2), although they only describe it in the case of $G$ finite,
being $\mathcal{F}$ the family of finite subgroups of $G$, and
with another purpose. In the language of that paper, we have
proved that $X$ is the $\mathcal{F}$-approximation to a point.

\end{rem}

We conclude this section with a modification of the previous
models that will be useful in the sequel. So, if $\{H_1,\ldots
,H_n\}$ are subgroups of $G$ that belong to $\mathcal{F}$ and such
that for every $0<i<n$ there exists a $G$-equivariant map
$G/H_i\longrightarrow G/H_{i+1}$, we define
$\mathbf{G/H_1\longrightarrow\ldots\longrightarrow G/H_n}$ as the
small category whose elements are $n$-uples $(a_1H_1,\ldots
,a_nH_n)$ such that for every $0<i<n$ there exists a
$G$-equivariant map $f_i:G/H_i\longrightarrow G/H_{i+1}$ with
$f_i(a_iH_i)=a_{i+1}H_{i+1}$, and whose morphisms are the identity
maps. Now, if $\mathbf{\Gamma}$ is the poset category of
non-degenerate simplices of
$\textrm{N}(\mathbf{O}_{\mathcal{F}})$, we define a functor
$S:\mathbf{\Gamma}^{op}\longrightarrow\mathbf{Cat}$ that takes the
simplex represented by the chain of maps
$\{G/H_1\longrightarrow\ldots\longrightarrow G/H_n\}$ to the
category $\mathbf{G/H_1\longrightarrow\ldots\longrightarrow
G/H_n}$, and the face maps to the obvious functors. Then we have
the following:

\begin{prop}

In the previous conditions,
$\emph{hocolim}_{\mathbf{\Gamma}}\emph{N}(S)$ is a model for
$\emph{E}_{\mathcal{F}}G$, and $\emph{N}(\mathbf{\Gamma})$ is a
model for $\emph{B}_{\mathcal{F}}G$.

\label{modelomodific}
\end{prop}

\begin{proof}

By theorem \ref{Thomason},
$\textrm{hocolim}_{\mathbf{\Gamma}}\textrm{N}(S)\simeq
\textrm{N}(\mathbf{Gr}(S))$, and on the other hand, the left
action of $G$ over every homogeneous space $G/H$,
($H\in\mathcal{F}$) induces, via functors, another one over the
categories $\mathbf{G/H_1\longrightarrow\ldots\longrightarrow
G/H_n}$. Now, observe that $\textrm{N}(\mathbf{Gr }(S))$ is the
subdivision (in the sense of \cite{Goerss99}, III.4) of
$\textrm{N}(\mathbf{Gr }(R))$, where $R$ is the functor defined in
proposition \ref{modelo}. In fact, if we take a non-degenerate
simplex of $\textrm{N}(\mathbf{Gr }(R))$ that is represented by a
chain of morphisms
$(G/H_1,a_1H_1)\longrightarrow\ldots\longrightarrow
(G/H_n,a_nH_n)$, its barycenter is the vertex of the nerve
$\textrm{N}(\mathbf{Gr }(S))$ represented by the object
$(\mathbf{G/H_1\longrightarrow\ldots\longrightarrow G/H_n},
a_1H_1$ $\longrightarrow\ldots\longrightarrow a_nH_n)$ of
$\textbf{Gr}(S)$. Then, by (\cite{Jardine02} prop. 12-14), there
exists a homotopy equivalence $|\textrm{N}(\mathbf{Gr }R)|\simeq
|\textrm{N}(\textbf{Gr }S)|$, that is a $G$-equivalence by
construction, and then $\textrm{N}(\mathbf{Gr }S)$ is a model for
$\textrm{E}_{\mathcal{F}}G$. A similar line of reason proves that
$\textrm{N}(S)$ is a model for $\textrm{B}_{\mathcal{F}}G$.
\end{proof}

The main advantage of these models is that they reconstruct
$\textrm{E}_{\mathcal{F}}G$ and $\textrm{B}_{\mathcal{F}}G$ as
homotopy colimits over a \emph{poset} category, and in particular
they have structure of simplicial \emph{complexes}. These facts
will be very useful in the final section.

\section{$\textrm{\underline{B}}G$ is homotopy equivalent to a nullification}

Let $G$ be a discrete group, $\mathcal{F}$ a family of subgroups
of $G$ closed under conjugation and taking subgroups. As a first
step in our study of the relation between the classical and proper
classifying spaces, we will describe a canonical map that always
relates $\textrm{B}G$ and $\textrm{B}_{\mathcal{F}}G$.

So, consider models of E$G$ and $\textrm{E}_{\mathcal{F}}G$; both
of them are $G$-spaces, and then we can make the Borel
construction $\textrm{E}G\times_{G}\textrm{E}_{\mathcal{F}}G$.
Now, let $p_1$ be the projection
$$\textrm{E}G\times_{G}\textrm{E}_{\mathcal{F}}G\stackrel{p_1}{\longrightarrow}
\textrm{E}G/G\simeq\textrm{B}G.$$ The action of $G$ over
$\textrm{E}G$ is free, so the map $p_1$ is a fibration, and its
homotopy fiber $\textrm{E}_{\mathcal{F}}G$ is contractible. Thus,
$p_1$ is a homotopy equivalence, and
$\textrm{E}G\times_{G}\textrm{E}_{\mathcal{F}}G$ is a model of
B$G$.

Consider now the projection over the second component
$$\textrm{E}G\times_{G}\textrm{E}_{\mathcal{F}}G\stackrel{p_2}{\longrightarrow}
\textrm{E}_{\mathcal{F}}G/G\simeq\textrm{B}_{\mathcal{F}}G.$$  We
have seen that $\textrm{E}G\times_{G}\textrm{E}_{\mathcal{F}}G$
and $\textrm{E}_{\mathcal{F}}G/G$ are respectively models for B$G$
and $\textrm{B}_{\mathcal{F}}G$, and then $p_2$ can be thought as
a map B$G\longrightarrow\textrm{B}_{\mathcal{F}}G$ that we will
call $f$ in the rest of the section. The map $f$ is not a
fibration in general, because the action of $G$ over
$\textrm{E}_{\mathcal{F}}G$ is not free. In fact, if
$x\in\textrm{B}_{\mathcal{F}}G$, we have that $f^{-1}(x)$ has the
homotopy type of $\textrm{E}G\times_{G}G/H_{x}$, being $H_{x}$ the
isotropy group of $x$, that belongs to $\mathcal{F}$. Now,
$\textrm{E}G\times_{G}G/H_{x}$ is a model for $\textrm{B}H_{x}$,
and hence all the fibers of the map $f$ have the homotopy type of
classifying spaces of groups of $\mathcal{F}$.

This fact gave us the intuition that the map we have studied could
encode a functorial way of passing from the usual classifying
space of $G$ to the classifying space for proper $G$-bundles, and
what is more important, to obtain valuable information of the
latter starting from the homotopy structure of $\textrm{B}G$, and
viceversa. More concretely, we searched for a functor $F$ such
that the following conditions hold:

\begin{enumerate}

\item $F$ ``kills" the homotopy fiber of $f$.

\item $F(f)$ is a weak equivalence.

\item $F(\underline{\textrm{B}}G)\simeq\underline{\textrm{B}}G.$

\end{enumerate}

The two first conditions give the impression, according to
(\cite{Dror-Farjoun95}, 1.H.1 and 3.D.3), of $F$ being a
localization functor $L$ in the sense of Dror-Farjoun, and in
fact, the functor we have found has been the $A$-nullification
functor with respect to a certain space $A$. For the main
properties of these functors you can look at \cite{Bousfield94},
\cite{Chacholski96} or \cite{Dror-Farjoun95}, although we will
recall here the definition.

Let $A$ and $X$ be spaces; $X$ is said $A$-\emph{null} if the
mapping space $\textrm{map}(A,X)$ is homotopy equivalent to $X$
via the inclusion of constant maps
$X\longrightarrow\textrm{map}(A,X)$. The $A$-nullification of $X$
is a functor
$\mathbf{P}_A:\mathbf{Spaces}\longrightarrow\mathbf{Spaces}$ that
takes every space $X$ to an $A$-null space $\mathbf{P}_AX$ such
that there exists a universal map $X\longrightarrow\mathbf{P}_AX$
which induces a weak homotopy equivalence
$$\textrm{map}(\mathbf{P}_AX,Y)\simeq\textrm{map}(X,Y)$$ for every $A$-null
space $Y$. It can be seen that any other $A$-null space $X$ for
which the last property holds is indeed homotopy equivalent to
$\mathbf{P}_AX$. So, we have defined a functor in the category of
unpointed spaces (that can also be defined in the category of
pointed spaces, although we we will only work here with unpointed
spaces), which is always coaugmented and idempotent, and kills the
structure of $X$ that ``depends'' on $A$. In fact, $\mathbf{P}_AX$
is the localization of $X$ with respect to the constant map
$A\longrightarrow {*}$, and the spaces $X$ for which
$\mathbf{P}_AX$ is contractible are called $A$-acyclic.

Consider now the set of all prime numbers $\{p_1,p_2,p_3\ldots\}$
with the usual order, and let $X$ be a space. In the remaining of
the paper we will denote by $W_n$ the space
$\textrm{B}\mathbb{Z}/p_1\vee\ldots\vee\textrm{B}\mathbb{Z}/p_n$,
and by $W_{\infty}$ the wedge $\bigvee\textrm{B}\mathbb{Z}/p$
extended over of the all prime numbers. The next key lemma is one
strong reason that suggests that the $W_{\infty}$-nullification is
the functor we need.

\begin{thm}

If $G$ is a finite group, then $\mathbf{P}_{W_{\infty}}\emph{B}G$
is contractible.

\end{thm}

\begin{proof}
A point is always null, so we only need to prove that  for every
$W_{\infty}$-null space $X$ there is an equivalence
$X\simeq\textrm{map}(\mathbf{P}_{W_{\infty}}\textrm{B}G,X)$. But
the $W_{\infty}$-null spaces are, in particular, $W_n$-null for
every $n$; hence,
$$\textrm{map}(\mathbf{P}_{W_{\infty}}\textrm{B}G,X)\simeq\textrm{map}(\textrm{B}G,X)\simeq
\textrm{map}(\mathbf{P}_{W_n}\textrm{B}G,X)$$ for every $n$. Now,
suppose that $|G|=p_{j_1}^{n_1}p_{j_2}^{n_2}\ldots p_{j_m}^{n_m}$,
with $j_1<\ldots j_m$. Because of (\cite{Yo}, 3.3), we know that
$\mathbf{P}_{W_k}\textrm{B}G$ is contractible for every $k\geq
j_m$. This implies
$$\textrm{map}(\mathbf{P}_{W_{\infty}}\textrm{B}G,X)\simeq
\textrm{map}(\mathbf{P}_{W_n}\textrm{B}G,X)\simeq X$$ as we
claimed.
 \end{proof}

Now suppose that $G$ is a discrete group; we are in position of
stating our main theorem:

\begin{thm}
Let $G$ be a group such that there exists a finite-dimensional
model for $\emph{\underline{B}}G$, and let $W_{\infty}$ be as
before; we have that $\emph{\underline{B}}G$ is homotopy
equivalent to $\mathbf{P}_{W_{\infty}}\emph{B}G.$

\label{bbarnull}

\end{thm}

\begin{proof}

Consider the model of $\textrm{\underline{E}}G$ given in section 2
as the realization of the nerve of the Grothendieck construction
of a functor
$$R:\mathbf{O}_{\mathcal{F}}\longrightarrow\mathbf{Cat}$$
described there. We have seen that
$\textrm{E}G\times_G\underline{\textrm{E}}G$ is a model of
$\textrm{B}G$, so using theorem \ref{Thomason} we obtain
$$\textrm{B}G\simeq\textrm{E}G\times_G\underline{\textrm{E}}G\simeq\textrm{E}G\times_G
(\textrm{hocolim}_{\mathbf{O}_{\mathcal{F}}}R)\simeq\textrm{hocolim}_{\mathbf{O}_{\mathcal{F}}}(\textrm{E}G\times_G
R(-))$$ where the previous equivalence is a simple application of
(\cite{Dwyer00}, 6.5). Now, if we apply the nullification functor
$\mathbf{P}_{W_{\infty}}$ to the previous string of equivalences,
we obtain a weak equivalence
$\mathbf{P}_{W_{\infty}}\textrm{B}G\simeq\mathbf{P}_{W_{\infty}}\textrm{hocolim}_{\mathbf{O}_{\mathcal{F}}}
(\textrm{E}G\times_G R(-))$, and the latter is equivalent, by
(\cite{Dror-Farjoun95} 1.H.1), to
$\mathbf{P}_{W_{\infty}}\textrm{hocolim}_{\mathbf{O}_{\mathcal{F}}}\mathbf{P}_{W_{\infty}}(\textrm{E}G\times_G
R(-))$. Observe that the spaces that appear in the target of the
functor $$\textrm{E}G\times_G
\textrm{N}(R(-)):\mathbf{O}_{\mathcal{F}}\longrightarrow\mathbf{Spaces}$$
have the homotopy type of classifying spaces of finite subgroups
of $G$. Hence, if we apply the previous proposition, we have that
$\mathbf{P}_{W_{\infty}}\circ (\textrm{E}G\times_G R(-))$ is
equivalent to the constant functor, and then
$$\mathbf{P}_{W_{\infty}}\textrm{B}G\simeq\mathbf{P}_{W_{\infty}}\textrm{hocolim}_{\mathbf{O}_{\mathcal{F}}}*
\simeq\mathbf{P}_{W_{\infty}}(\textrm{N}(\mathbf{O}_{\mathcal{F}}))\simeq\mathbf{P}_{W_{\infty}}(\underline{\textrm{B}}G).$$
Then, by the solution of Miller to the Sullivan conjecture
\cite{Miller84}, we know that the space
$\textrm{map}(W_{\infty},\textrm{\underline{B}}G)$ is homotopy
equivalent to $\textrm{\underline{B}}G$, and hence
$\textrm{\underline{B}}G$ is $W_{\infty}$-null. This means that
$\mathbf{P}_{W_{\infty}}\textrm{\underline{B}}G$ is homotopy
equivalent to $\textrm{\underline{B}}G$, and we are done.
\end{proof}

The following generalization, that will be widely applied in
section 6, is an immediate consequence of the proof of the
previous theorem:

\begin{cor}

If $\mathcal{F}$ is a family of finite subgroups of $G$ closed under conjugation and taking
 subgroups, and $\mathbf{P}_A$ is a nullification functor such that
 $\mathbf{P}_A\emph{B}H\simeq *$ for every $H\in\mathcal{F}$, then the map
 $f:\emph{B}G\longrightarrow\emph{B}_{\mathcal{F}}G$ is an equivalence after
 $A$-nullification.

\end{cor}

If we particularize for the family of all the finite groups, we obtain the following:

\begin{cor}

If $G$ is a discrete group, the classifying spaces $\emph{B}G$ and
$\underline{\emph{B}}G$ are always equivalent after
$W_{\infty}$-nullification; moreover, the map $f$ that was
described at the beginning of this section is, in fact, equivalent
to the $W_{\infty}$-nullification map if $G$ admits a
finite-dimensional model for $\underline{\emph{B}}G$.
\label{BGyBpG}
\end{cor}

\begin{proof}

The fact that $\textrm{B}G$ and $\textrm{B}_{\mathcal{F}}G$ are
$\mathbf{P}_{W_{\infty}}$-equivalent is a particular case of the
previous corollary. For the second statement, if we have the
fibration $$\mathbf{Fib
}f\longrightarrow\textrm{B}G\longrightarrow
\textrm{B}_{\mathcal{F}}G$$ (where $\mathbf{Fib} f$ stand for the
homotopy fibre of $f$), the base is $W_{\infty}$-null, and then by
(\cite{Dror-Farjoun95}, 3.D.3) the fibration is preserved under
$W_{\infty}$-nullification. Now the result is an easy consequence
of (\cite{Chacholski96}, 14.2).  \end{proof}

It is worth to point out that the finiteness conditions under
which the main theorem holds is that there exists a model of
$\textrm{\underline{B}}G$ for which $\textrm{map}_*(W_{\infty},
\textrm{\underline{B}}G)$ is weakly contractible. This is weaker
than having a finite-dimensional model for
$\textrm{\underline{B}}G$, but the condition that will always hold
for the groups that appear in the rest of this note will be the
latter, because it is the usual one that appears in the literature
of Geometric Group Theory. It would be interesting to find
cohomological conditions to be a Miller space (that is, spaces $X$
for which $\textrm{map}(\textrm{B}\mathbb{Z}/p,X)\simeq X)$ for
all the primes at the same time, because we could then apply the
theorem to these spaces. Recall that there already exist
well-known cohomological conditions of this kind for isolated
primes, as for example the Lannes-Schwartz Theorem
\cite{Lannes89}.

On the other hand, we cannot expect that theorem \ref{bbarnull}
holds for any $G$. To see this, take a space $X$ which is not
$W_{\infty}$-null. According to \cite{Leary01}, there exists a
discrete group $G$ such that $X$ is a model for
$\textrm{\underline{B}}G$; so our result does not hold for $G$.

We finish this section by showing that the discreteness of $G$ is
necessary in theorem \ref{bbarnull}.

\begin{exmp}

Let us consider the classifying space of $S^1$. As the circle is
compact, the classifying space for proper $S^1$-bundles is a point
by definition. On the other hand, consider the rationalization map
$\emph{B}S^1\longrightarrow K(\mathbb{Q},2)$. By the homotopy long
exact sequence, the homotopy fibre of this map has the homotopy
type of $ K(\oplus\mathbb{Z}_{p^{\infty}},1)$, where the direct
sum runs over all primes. As every Pr{\"u}fer group
$\mathbb{Z}_{p^{\infty}}$ is a colimit of a telescope of
injections between $p$-groups, the results (\cite{Dror-Farjoun95}
1.D.3) and (\cite{Yo}, 3.3) imply that
$\mathbf{P}_{W_{\infty}}K(\oplus\mathbb{Z}_{p^{\infty}},1)$ is
contractible. Then, by (\cite{Dror-Farjoun95} 1.H.1) the
rationalization is preserved by $W_{\infty}$-nullification, and as
$K(\mathbb{Q},2)$ is clearly $W_{\infty}$-null, we have that
$\mathbf{P}_{W_{\infty}}\emph{B}S^1\simeq K(\mathbb{Q},2)$, which
is non-contractible. In fact, it seems plausible to conjecture
that the $W_{\infty}$-nullification of the classifying space of a
compact Lie group is homotopy equivalent to its rationalization.

\end{exmp}

\section{The homotopy type of $\textrm{\underline{B}}G$}

In this paragraph we are going to prove some interesting
consequences that the main theorem \ref{bbarnull} has over the
homotopy type of $\textrm{\underline{B}}G$. Essentially, the idea
is to use properties of the nullification functors for describing
the classifying space for $G$-proper bundles.

\begin{rem}

From now on, we particularize for the case $\mathcal{F}$ being the
family of \emph{all} the finite subgroups of $G$, although a great
part of the results we obtain in the next sections remains valid
for any subfamily of $\mathcal{F}$ that is subgroup-closed and
conjugation-closed. Even if we do not mention it, we will also
suppose that the finiteness conditions of theorem \ref{bbarnull}
hold for all the groups that appear in this section.

\end{rem}

We begin by analyzing the behaviour of the functor
$\underline{\textrm{B}}$ under products.

\begin{prop}

Let $G_1$ and $G_2$ discrete groups which possess a finite
dimensional model for the classifying space for proper bundles.
Then the following holds:

\begin{itemize}

\item A model for $\emph{\underline{B}}(G_1\times G_2)$ is given
by $\emph{\underline{B}}G_1\times\emph{\underline{B}}G_2$.

\item The wedge
$\emph{\underline{B}}G_1\vee\emph{\underline{B}}G_2$ is a model
for $\emph{\underline{B}}(G_1{*}G_2)$.

\end{itemize}

\end{prop}

\begin{proof}

It is known that $\textrm{B}(G_1\times
G_2)\simeq\textrm{B}G_1\times\textrm{B}G_2$. Using that
$\textrm{\underline{B}}G_1\times\textrm{\underline{B}}G_2$ is
$W_{\infty}$-null (because the finiteness) and the preservation
property (\cite{Dror-Farjoun95} 1.A.8, prop. 4), we obtain that
$$\textrm{\underline{B}}(G_1\times G_2)\simeq\mathbf{P}_{W_{\infty}}(\textrm{\underline{B}}
(G_1\times
G_2))\simeq\mathbf{P}_{W_{\infty}}(\textrm{\underline{B}}G_1)\times\mathbf{P}_{W_{\infty}}(\textrm{\underline{B}}G_2)
\simeq\textrm{\underline{B}}G_1\times\textrm{\underline{B}}G_2.$$
The proof of the second statement is similar, using that
$\textrm{B}(G_1{*}G_2)\simeq\textrm{B}G_1\vee\textrm{B}G_2$ and
recalling that we can apply (\cite{Dror-Farjoun95}, 1.D.5) because
a wedge of a special case of pointed homotopy colimit.
\end{proof}

It is worth to point out that these results are established in
(\cite{Lueck04}, 4.9) without any assumption on the dimension,
although this way of proving them is probably new. On the other
hand, the second of them can be generalized to other colimits,
like some telescopes and pushouts.

\begin{prop}

Let $\{G_i\}_{i\in\mathbb{N}}$ be a collection of discrete groups
which possess a finite-dimensional model for the classifying space
for proper bundles . Then the following holds:

\begin{itemize}

\item If we have the pushout of groups $$ \xymatrix{ G_1
\ar[r]^{\alpha} \ar[d]_{\beta} & G_2 \ar[d] \\ G_3 \ar[r] & G }
$$ and the homomorphisms $\alpha$ and $\beta$ are injective, then
the pushout of the induced diagram of classifying spaces for
proper $G$-bundles is a model for $\underline{\emph{B}}G$.

\item If $G_1\longrightarrow G_2\longrightarrow
G_3\longrightarrow\ldots$ is a telescope of groups where the maps
are injective and we denote by $G$ the colimit of the telescope,
we have that the colimit of the telescope induced by
$\underline{\emph{B}}$ is a model for $\emph{\underline{B}}G$.

\end{itemize}

\end{prop}

\begin{proof}

To prove the first statement, recall that by Whitehead's theorem
(\cite{Brown82}, II.7.3) the pushout of the classical classifying
spaces is the classifying space of the pushout. As the inclusions
$\textrm{B}G_1\hookrightarrow\textrm{B}G_2$ and
$\textrm{B}G_1\hookrightarrow\textrm{B}G_3$ are cofibrations,
$\textrm{B}G$ has the homotopy type of the homotopy pushout. If we
apply now the functor $\mathbf{P}_{W_{\infty}}$ to the diagram,
the result is deduced from theorem \ref{bbarnull},
(\cite{Dror-Farjoun95}, 1.D.3) and the fact that there exists a
finite-dimensional model for the homotopy pushout of the induced
diagram
$$\textrm{\underline{B}}G_2\longleftarrow\textrm{\underline{B}}G_1\longrightarrow\textrm{\underline{B}}G_3.$$
The second statement can be proved in an analogous way using again
the relationship between localization and colimits given in
(\cite{Dror-Farjoun95}, 1.D.3) and the fact that the strict
colimit of a telescope of cofibrations has always the homotopy
type of the homotopy colimit.  \end{proof}

As Whitehead's theorem is not true if the maps that appear in the
diagram are not injective, it should not be expected that the
functor $\textrm{\underline{B}}$ preserve colimits in full
generality.

Recall now that for any discrete group $G$ the fundamental group
of $\textrm{\underline{B}}G$ can be identified (\cite{Leary01},
prop. 3) as the quotient of $G$ by the (normal) subgroup generated
by the torsion elements. Using the main theorem \ref{bbarnull}, we
can identify in some cases the universal cover of
$\textrm{\underline{B}}G$.

\begin{prop}

Let $G$ be a discrete group which has a finite-dimensional model
for $\emph{\underline{B}}G$, and let $T<G$ be the subgroup
generated by the torsion elements. If the quotient $G/T$ is
torsion-free, then the universal cover of $\emph{\underline{B}}G$
has the homotopy type of the $W_{\infty}$-nullification of
$\emph{B}T$.

\end{prop}

\begin{proof}

We know that $T$ is normal in $G$, so we have a fibration
$$\textrm{B}T\longrightarrow\textrm{B}G\longrightarrow\textrm{B}(G/T).$$
As $G/T$ is torsion-free, its classifying space is
$W_{\infty}$-null. Thus the previous fibration is preserved by
$W_{\infty}$-nullification, and we obtain another one:
$$\mathbf{P}_{W_{\infty}}\textrm{B}T\longrightarrow\textrm{\underline{B}}G\longrightarrow\textrm{B}(G/T).$$
Note that, as $T$ is a subgroup of $G$, every model for
$\textrm{\underline{E}}G$ is also a model for
$\textrm{\underline{E}}T$. Hence, $\textrm{\underline{B}}T$ is a
model for $\mathbf{P}_{W_{\infty}}\textrm{B}T$, and in particular
$\pi_1(\mathbf{P}_{W_{\infty}}\textrm{B}T)=\pi_1(\textrm{\underline{B}}T)=\{1\}$.
This implies that $\mathbf{P}_{W_{\infty}}\textrm{B}T$ is
simply-connected, and we are done.  \end{proof}

The last important consequence of the main theorem that we are
going to prove here has to do with fibrations, and will have great
importance in the remaining of this note.

It is a well-known fact of basic algebraic topology that if we
have a group extension, then the sequence induced at the level of
classifying spaces is a fibration sequence. Using the description
of theorem \ref{bbarnull}, we find sufficient conditions that
guarantee that the analogous result for $\underline{\textrm{B}}G$
holds, and we show by means of an easy example that the statement
does not need to be true if those hypotheses fail to be fulfilled.

So, suppose we have a short exact sequence of groups
$$\{1\}\longrightarrow G_1\longrightarrow G\longrightarrow
G_2\longrightarrow\{1\}$$ which have a finite-dimensional model
for the classifying space for proper bundles. Then the following
result is true:

\begin{prop}

If $G_2$ is torsion-free or $G_1$ admits a contractible model for
$\underline{\emph{B}}G_1$, the homotopy fiber of the induced map
$\underline{\emph{B}}G\longrightarrow\underline{\emph{B}}G_2$ is
homotopy equivalent to $\underline{\emph{B}}G_1.$

\end{prop}

\begin{proof}

It is enough to combine the results (\cite{Dror-Farjoun95}, 1.H.1
and 3.D.3) of Dror-Farjoun with our description of
$\underline{\textrm{B}}G$ as a nullification (Theorem
\ref{bbarnull}).  \end{proof}

The next example, that was one of the first motivations for our
work, will show that the above conditions are necessary.

\begin{exmp}
\label{toro}

Consider the group $D_{\infty}\times D_{\infty}$, and $H$ the index
two subgroup whose elements are the words that can be written with
an even number of letters. We have an extension
$$H\longrightarrow D_{\infty}\times
D_{\infty}\longrightarrow\mathbb{Z}/2$$ that induces a sequence of
maps
$$\underline{\textrm{B}}H\longrightarrow\underline{\textrm{B}}(D_{\infty}\times
D_{\infty})\longrightarrow\underline{\textrm{B}}\mathbb{Z}/2.$$ It
is not hard to see that a model for
$\underline{\textrm{E}}(D_{\infty}\times D_{\infty})$ is given by
$\mathbb{R}^2$, and the quotient by the action of
$(D_{\infty}\times D_{\infty})$ is a square, which is
contractible. By (\cite{Leary01}, prop. 8, see also lemma
\ref{lemawall} below), $\underline{\textrm{B}}H$ is homotopy
equivalent to the 2-sphere, and on the other hand, $\mathbb{Z}/2$
is finite and so $\underline{\textrm{B}}\mathbb{Z}/2$ is
contractible. This means that the aforementioned sequence cannot
be a fibration sequence. Of course, neither of the conditions of
the proposition hold in this case.

Observe that $D_{\infty}\times D_{\infty}$ is a group of
isometries of the plane. We will carefully study these groups in
section 6.

\label{dihinf}

\end{exmp}

\section{Homotopy models of $\textrm{\underline{B}}G$ for some classes of discrete groups}
\label{discrete}

In this section we will use the theorem \ref{bbarnull} for
describing the homotopy type of $\underline{\textrm{B}}G$ for a
wide range of groups. As a by-product, we will obtain that for
every $G$ nilpotent such that it admits a finite-dimensional model
for $\underline{\textrm{B}}G$, $\underline{\textrm{B}}G$ is
nilpotent as a space, and we also determine, for $p$ odd, the
$\textrm{B}\mathbb{Z}/p$-nullification of classifying spaces of
supersoluble groups. Let us start by considering the class of
locally finite groups.

\subsection{Locally finite groups}
\label{locafini}

It is known that the classifying space for proper $G$-bundles of a
group $G$ is contractible if the group $G$ is locally finite. We
begin this section by presenting an easy proof of this fact in an
ample range of cases.

\begin{prop}

Let $G$ be a locally finite group that admits a finite-dimen\-sional
model for $\underline{\emph{B}}G$. Then $\underline{\emph{B}}G$ is
contractible. \label{locally}

\end{prop}

\begin{proof}

It is known (see for example \cite{Miller84}, 9.8) that every
locally finite group is the colimit of the directed system of its
finite subgroups. Thus, we have a homotopy equivalence
$\textrm{B}G\simeq\textrm{hocolim}_{\mathcal{C}}\textrm{B}H$,
where $\mathcal{C}$ is a contractible poset category (because it
has an initial object given by the trivial group) and
$\textrm{B}H$ represents all the classifying spaces of finite
groups $H$ of $G$. So, by (\cite{Dror-Farjoun95}, 1.D.3), we
obtain
$$\mathbf{P}_{W_{\infty}}\textrm{B}G\simeq\mathbf{P}_{W_{\infty}}(\textrm{hocolim}_{\mathcal{C}}
\mathbf{P}_{W_{\infty}}\textrm{B}H)=\mathbf{P}_{W_{\infty}}(|\mathcal{C}|)=\mathbf{P}_{W_{\infty}}(*)=*$$
and we are done. \end{proof}

\begin{rem}

By \cite{DKLT02}, this result applies to all locally finite groups
whose cardinal is smaller than $\aleph_{\omega}$.
\label{locfinite}
\end{rem}

Now we can prove the following result, that concerns to the
classifying space for proper $G$-bundles of extensions of locally
finite groups.

\begin{prop}

Let $$\{1\}\longrightarrow K\longrightarrow G\longrightarrow
Q\longrightarrow\{1\}$$ be an extension of groups, $K$ a locally
finite group whose cardinal is smaller than $\aleph_{\omega}$, and
assume there is a bound on the order of finite subgroups of $Q$.
Then if $Q$ admits a finite model for $\underline{\emph{B}}Q$, $G$
admits a finite model for $\underline{\emph{B}}G$, and then
$\underline{\emph{B}}Q$ is homotopy equivalent to
$\underline{\emph{B}}G$.

\end{prop}

\begin{proof}

If we apply the results \ref{bbarnull} and \ref{locally} we obtain
the statement is true if there is a finite-dimensional model for
$\underline{\textrm{B}}G$, and this happens by proposition 4.4 of
\cite{Mislin01}.
\end{proof}

\subsection{Groups with the normalizer condition}

We study now the groups for which the normalizer condition holds.
It is greatly remarkable that this class contains all the
nilpotent groups.

Recall that a group $G$ is said to satisfy the \emph{normalizer
condition} if every proper subgroup of $G$ is distinct from its
normalizer. In this case the following holds (see \cite{Kurosh60},
page 215):

\begin{enumerate}

\item For every prime $p$, there exists a normal $p$-group $T_p$
such that if $x\in G$ and the order of $x$ is a power of $p$, then
$x\in T_p$.

\item The elements of finite order of $G$ form a normal subgroup
of $G$ which is isomorphic to $\prod_{p\textrm{ }prime}T_p$.

\end{enumerate}

Throughout this section we will impose to the groups for which the
normalizer condition holds that the torsion $p$-subgroups $T_p$
that we have just defined are locally finite. We need this
condition because it is not known if the
$\textrm{B}\mathbb{Z}/p$-nullification of the classifying space of
a general $p$-group is contractible if the group is not locally
finite. Among the few examples that have been described of
non-locally finite $p$-groups we can remark the Burnside groups
$B(n,e)$ for $n>1$ and $e>664$ or the ``monsters" of
Tarski-Olshanskii. See \cite{Adyan75} and \cite{Olshanskii89} for
more information about these families of groups.

All these facts have the following interesting consequence:

\begin{prop}

If $G$ is a discrete group for which the normalizer condition
holds, and $p_1\ldots p_n$ is a collection of primes, we have a
homotopy equivalence $\mathbf{P}_{B\mathbb{Z}/p_1\vee\ldots\vee
B\mathbb{Z}/p_n}\emph{B}G\simeq\emph{B}(G/T_{p_1}\times\ldots
\times T_{p_n}).$ In particular, there is an equivalence
$\mathbf{P}_{W_{\infty}}\emph{B}G\simeq\emph{B}(G/\prod_{p\textrm{ }prime}T_{p}).$ 

\label{normaliser}
\end{prop}

\begin{proof}

For simplicity, we will only prove the case of one prime $p$ (the
generalization to a family is immediate). It is clear that
$\textrm{B}T_p$ is $\textrm{B}\mathbb{Z}/p$-acyclic, and
$\textrm{B}(G/T_p)$ is $\textrm{B}\mathbb{Z}/p$-null, so if we
$\textrm{B}\mathbb{Z}/p$-nullify the fibration
$$\textrm{B}T_p\longrightarrow\textrm{B}G\longrightarrow\textrm{B}(G/T_p)$$
we obtain the desired homotopy equivalence.  \end{proof}

If we suppose that $G$ is such that exists a finite-dimensional
model for $\underline{\textrm{B}}G$ we have:

\begin{cor}

$\underline{\emph{B}}G\simeq\emph{B}(G/\prod_{p\textrm{
}prime}T_{p}).$

\end{cor}

So we have a complete description of the homotopy type of
$\underline{\textrm{B}}G$. \newline

Other case that can be solved with the same tools is the
following:

\begin{prop}
Let $G$ be a discrete group, $H$ a normal subgroup of $G$ for
which the normalizer conditions holds, and such that $G/H$ does
not have $p$-torsion. If $T_p$ is the $p$-torsion subgroup of $H$,
then the $\emph{B}\mathbb{Z}/p$-nullification of $\emph{B}G$ fits
into the following fibration sequence:
$$\emph{B}(H/T_p)\longrightarrow\mathbf{P}_{B\mathbb{Z}/p}\emph{B}G\longrightarrow\emph{B}(G/H)$$
and hence it is an Eilenberg-MacLane space.

\end{prop}

\begin{proof}

The base of the fibration
$$\textrm{B}H\longrightarrow\textrm{B}G\longrightarrow\textrm{B}(G/H)$$
is $\textrm{B}\mathbb{Z}/p$-null, so by (\cite{Dror-Farjoun95},
3.D.3) the fibre is preserved under
$\textrm{B}\mathbb{Z}/p$-nullification. The result now follows
from proposition \ref{normaliser}.  \end{proof}

Taking into account the main theorem \ref{bbarnull}, the following
corollary is immediate:

\begin{cor}

In the hypotheses of the previous proposition, if $G/H$ is
tor\-sion-free, $T$ is the torsion subgroup of $H$ and there exists
a finite-dimensional model for $\emph{\underline{B}}G$ , then the
fibration
$$\emph{B}(H/T)\longrightarrow\underline{\emph{B}}G\longrightarrow\emph{B}(G/H)$$
defines the classifying space for proper $G$-bundles, which is
again an Eilenberg-MacLane space.

\end{cor}

We conclude this paragraph by focusing on nilpotent groups, that
is a distinguished class of discrete groups for which the
normalizer condition holds. The following result proves that the
$\textrm{B}\mathbb{Z}/p$-nullification preserves nilpotency when
it is applied on classifying spaces of nilpotent groups, and in
fact, the functor $\underline{\textrm{B}}$ sends nilpotent groups
(for which the finiteness conditions hold) to nilpotent spaces.

\begin{cor}

If $G$ is a nilpotent group, the nullification
$\mathbf{P}_{B\mathbb{Z}/p_1\vee\ldots\vee
B\mathbb{Z}/p_n}\emph{B}G$ is, for every set of primes
$\{p_1,\ldots ,p_n\}$, the classifying space of a nilpotent group.
If moreover $G$ admits a finite-dimensional model for
$\underline{\emph{B}}G$, we obtain that $\underline{\emph{B}}G$ is
again the classifying space of a nilpotent group, and hence
nilpotent as a space.

\end{cor}

\begin{proof}

Using the previous results, it is enough recalling that the
quotient of a nilpotent group is always nilpotent, and that
according to (\cite{Olshanskii89}, 2.7.1), every nilpotent
$p$-group is locally finite.  \end{proof}

In particular, using (\cite{Mislin01}, 4.5) we have that the part
of the previous corollary that alludes to
$\underline{\textrm{B}}G$ is always true if $G$ is a nilpotent
group whose cardinal is smaller than $\aleph_{\omega}$ and whose
torsion-free rank is finite.

\subsection{Supersoluble groups}

In this paragraph we will compute, for $p$ odd, the
$\textrm{B}\mathbb{Z}/p$-nullification of classifying spaces of
supersoluble groups. In this case we obtain no result about the
homotopy type of $\textrm{\underline{B}}G$ (for reasons that will
be explained at the end) but we include here this computation
because of its intrinsical interest, and because the way it has
been worked out generalizes in some sense the methods we have used
to compute the $\textrm{B}\mathbb{Z}/p$-nullification in the
previous sections.

Recall that a group $G$ is called \emph{supersoluble} if it has
cyclic normal series of finite length. It is known that every
finitely generated nilpotent group is supersoluble, and that every
supersoluble group is polycyclic.

Our key result for computing
$\mathbf{P}_{{B}\mathbb{Z}/p}\textrm{B}G$ is the following
(\cite{Robinson72}, page 67):

\begin{prop}

If $G$ is a supersoluble group, there exists a characteristic
series $1\unlhd L\unlhd M\unlhd G$, in such a way that $L$ is
finite with odd order, $M/L$ is a finitely-generated torsion-free
nilpotent group and $G/M$ is a finite 2-group.

\end{prop}

In the sequel we will use the notation of this proposition. Let
$p$ be an odd prime, and consider the fibration:
$$\textrm{B}L\longrightarrow\textrm{B}M\longrightarrow\textrm{B}(M/L).$$ As $M/L$ is torsion-free, its classifying space is automatically
$\textrm{B}\mathbb{Z}/p$-null, and then by (\cite{Dror-Farjoun95},
3.D.3) we have the nullified fibration:
$$\mathbf{P}_{{B}\mathbb{Z}/p}\textrm{B}L\longrightarrow\mathbf{P}_{{B}\mathbb{Z}/p}\textrm{B}M
\longrightarrow\textrm{B}(M/L).$$ Using (\cite{Yo}, 3.3), the
fundamental group of $\mathbf{P}_{{B}\mathbb{Z}/p}\textrm{B}M$ is
identified by an extension
$$L/T_{\mathbb{Z}/p}L\longrightarrow\pi_1\mathbf{P}_{{B}\mathbb{Z}/p}\textrm{B}M\longrightarrow M/L$$
where $T_{\mathbb{Z}/p}L$ is the minimal normal subgroup of $L$
that contains all the $p$-torsion (the $\mathbb{Z}/p$-radical),
and the universal cover of
$\mathbf{P}_{{B}\mathbb{Z}/p}\textrm{B}M$ is homotopy equivalent
to $\mathbb{Z}[1/p]_{\infty}(\textrm{B}T_{\mathbb{Z}/p}L)$, where
$\mathbb{Z}[1/p]_{\infty}$ denotes Bousfield-Kan
$\mathbb{Z}[1/p]$-completion (see \cite{Bousfield72} for a
definition).
\newline

Now, we have the fibration that involves $M$ and $G$:
$$\textrm{B}M\longrightarrow\textrm{B}G\longrightarrow\textrm{B}(G/M).$$
This fibration is again preserved under
$\textrm{B}\mathbb{Z}/p$-nullification, because $G/M$ is a 2-group
and $p$ is odd. The long exact sequence of the nullified fibration
proves that the fundamental group of
$\mathbf{P}_{{B}\mathbb{Z}/p}\textrm{B}G$ fits into the following
exact sequence:
$$\pi_1\mathbf{P}_{{B}\mathbb{Z}/p}\textrm{B}M\longrightarrow\pi_1\mathbf{P}_{{B}\mathbb{Z}/p}
\textrm{B}G\longrightarrow G/M$$ where the kernel has already been
described. On the other hand, the universal cover of
$\mathbf{P}_{{B}\mathbb{Z}/p}\textrm{B}G$ is the same as the
universal cover of $\mathbf{P}_{{B}\mathbb{Z}/p}\textrm{B}M$ which
is $\mathbb{Z}[1/p]_{\infty}\textrm{B}T_{\mathbb{Z}/p}L$, as we
said before. Thus, we have described the desired nullification by
means of a covering fibration. On the other hand, the fact that
the classifying space of the quotient $G/M$ is \emph{not}
$\textrm{B}\mathbb{Z}/2$-null makes these methods useless for
computing the $\textrm{B}\mathbb{Z}/2$-nullification of
$\textrm{B}G.$ As an easy consequence of this, $\textrm{B}(G/M)$
is not $W_{\infty}$-null (in the notation of theorem
\ref{bbarnull}) and then we cannot get any homotopical description
of $\underline{\textrm{B}}G$ in this way.

\section{Nullifying classifying spaces of groups of isometries via proper actions}

So far we have applied theorem \ref{bbarnull} for obtaining
results about $\underline{\textrm{B}}G$ using properties of the
nullification functors. In this section we will go the other way
round, using geometric characteristics of the group $G$ for
describing topological features of the classifying space.

Our analysis has been focused in some of the crystallographic
groups of the plane, also known as wallpaper groups. Recall that
these are groups of isometries of $\mathbb{R}^2$ that fix a
pattern that is invariant under translations in the directions of
two lineally independent vectors. It is known that they are
exactly seventeen of these groups, and they are always finite
extensions of $\mathbb{Z}\oplus\mathbb{Z}$ by a finite group. The
main references available about the structure of these groups are
\cite{Schatts78} (that has been specially interesting for us
because the big amount of pictures fundamental domains, mirror
lines, rotation centers, generating regions, etc. that can be
found on it), \cite{Levy}, \cite{Xahlee}, \cite{Coxeter65},
\cite{DuSautuoy99} and \cite{Conway90}; we refer the reader to
them for the details of the structure of the groups that in the
sequel will stand without any explicit proof.

The general idea is to describe, for a prime $p$ and a discrete
group $G$ that has $p$-torsion, the homotopy type of the
$\textrm{B}\mathbb{Z}/p$-nullification of the classifying space of
$G$ using the main theorem \ref{bbarnull}. We have chosen
wallpaper groups essentially for two reasons: the first of them is
the following structure result, that is a particular case of
proposition (\cite{Lueck00}, 5.2):

\begin{lem}

Let $G$ be one of the seventeen wallpaper groups. Then
$\mathbb{R}^2$, endowed with the natural action of $G$, is a model
for $\underline{\emph{E}}G$.

\label{lemawall}
\end{lem}

The second feature of the wallpaper groups that we are going to
use is that all of them possess a well-known model for the orbit
space $\mathbb{R}^2/G$, which in fact is always described as a
finite-dimensional orbifold. A list of these standard models can
be found in \cite{Levy}. According to the previous lemma, these
spaces can be also interpreted as models for
$\underline{\textrm{B}}G$, and using this we will apply theorem
\ref{bbarnull} for obtaining the value of the
$\textrm{B}\mathbb{Z}/p$-nullification of the classifying space of
$G$.

Now, if $G$ is a wallpaper group that only has torsion in a family
of primes $P=\{p_1\ldots p_r\}$, it is an easy consequence of
corollary \ref{BGyBpG} that every model for
$\underline{\textrm{B}}G$ is a model for
$\textrm{B}_{\mathcal{F}}G$, being $\mathcal{F}$ the family of
finite subgroups of $G$ whose order is divided only for primes of
the family $P$. In particular, if $G$ has torsion only in $P$ and
admits a finite-dimensional model for $\underline{\textrm{B}}G$,
we have $\underline{\textrm{B}}G=
\textrm{B}_{\mathcal{F}}G\simeq\mathbf{P}_{\textrm{B}\mathbb{Z}/{p_1}\vee\ldots\vee
\textrm{B}\mathbb{Z}/{p_r}}\textrm{B}G$. We will frequently use
this fact in the sequel.

We are going to study here three of the wallpaper groups, namely
$\mathbf{pmm}$, $\mathbf{p3}$ and $\mathbf{p3m1}$. The main reason
of our choice is that they give examples in which the
$\textrm{B}\mathbb{Z}/p$-nullification of the classifying space
has homological dimension zero, positive and infinite,
respectively.

\medskip

(1)\qua \textbf{The discope group pmm}\qua As a group of symmetries of the
real plane, this group is generated by two perpendicular
translations and two reflections whose axes are also
perpendicular, and in fact it is isomorphic to $D_{\infty}\times
D_{\infty}$. Recall that this group already appeared in the
construction of example \ref{toro}.

The group \textbf{pmm} contains reflections and rotations, and as
it can be seen in the tables of \cite{Schatts78}, it tessellates
the plane with rectangles. The orbit space of the plane by the
action of this group is also a rectangle, and hence the
classifying space for $\mathbf{pmm}$-proper bundles is
contractible. On the other hand, as the rotations that appear in
the group are of order two, the group only contains torsion at the
prime two, and thus $\underline{\textrm{B}}\mathbf{pmm}$ is a
model for $\textrm{B}_{\mathcal{F}_2}\textrm{pmm}$. Now, according
to the main theorem \ref{bbarnull}, we conclude that
$\mathbf{P}_{B\mathbb{Z}/2}\textrm{B}\mathbf{pmm}\simeq {*}$.

\medskip

(2)\qua \textbf{The tritrope group p3}\qua This group is generated by two
translations whose directions form an angle of $\pi/3$ and a
rotation of angle $2\pi/3$. A presentation with these generators
is given by:
$$\mathbf{p3}=\{x,y,z;xyx^{-1}y^{-1}=1, z^3=1, zxz^{-1}y^{-1}x=1, zyz^{-1}x=1\}.$$
The distinguished isometries of this group are the 3-rotations, so
we have no reflections nor glide-reflections and the torsion is
concentrated in the prime three. The fundamental region of
$\mathbf{p3}$ (that is, the smallest region of $\mathbb{R}^2$
whose images under the action of \textbf{p3} cover the plane) is a
rhombus, and the action gives rise to a tessellation of
$\mathbb{R}^2$ by hexagons; in fact, this is the simplest
wallpaper group such that the induced tessellation is not by
quadrilaterals. The quotient $\mathbb{R}^2/\mathbf{p3}$ has then
the shape of a non-slit turnover with three corners and no mirror
points (see \cite{Xahlee} for details), and in particular it has
the homotopy type of a 2-sphere. Hence, using an analogous
argument to that of the previous case, we obtain that the
$\textrm{B}\mathbb{Z}/3$-nullification of $\textrm{B}\mathbf{p3}$
is homotopy equivalent to $S^2$.

\medskip

(3)\qua \textbf{The tryscope group p3m1}\qua A convenient system of
generators for this group is given by the two usual translations
(whose directions form again an angle of $\pi/3$), a rotation of
angle $2\pi/3$, and a reflection whose axis is the bisectrix of
the vectors that determine the generating translations; in
particular, the reflection gives torsion in the prime two and the
rotation gives it in the prime three.

On the other hand, maybe the best way to understand this group is
as the group of reflections in the sides of an equilateral
triangle. Hence, this is the fundamental region,  and the lattice
is hexagonal, as in the previous example. As one can see in
\cite{Schatts78}, the orbit space by the action of $\mathbf{p3m1}$
on $\mathbb{R}^2$ is a triangle, and as this is a model for the
classifying space for proper bundles, we have that
$\textrm{\underline{B}}\mathbf{p3m1}$ is contractible. Applying
one more time theorem \ref{bbarnull}, we obtain that the
$\textrm{B}\mathbb{Z}/2\vee\textrm{B}\mathbb{Z}/3$-nullification
of the classifying space of $\mathbf{p3m1}$ is a point.

Now we are also interested in the
$\textrm{B}\mathbb{Z}/3$-nullification of
$\textrm{B}\mathbf{p3m1}$, and we need to use a different
strategy. The tryscope group can be seen as an extension of
$\mathbb{Z}\oplus\mathbb{Z}$ by the symmetric group $\Sigma_3$,
and a consequence of this is that $\mathbf{p3}$ is an index two
subgroup of $\mathbf{p3m1}$. In particular, this gives rise to a
fibration:
$$\textrm{B}\mathbf{p3}\longrightarrow\textrm{B}\mathbf{p3m1}\longrightarrow\textrm{B}\mathbb{Z}/2.$$
The base space is $\textrm{B}\mathbb{Z}/3$-null, so according to
(\cite{Dror-Farjoun95}, 3.D.3) and our previous description of
$\mathbf{P}_{B\mathbb{Z}/3}\textrm{B}\mathbf{p3}$, the
$\textrm{B}\mathbb{Z}/3$-nullification of
$\textrm{B}\mathbf{p3m1}$ is identified by a covering fibration:
$$S^2\longrightarrow\mathbf{P}_{\textrm{B}\mathbb{Z}/3}\textrm{B}\mathbf{p3m1}\longrightarrow\textrm{B}\mathbb{Z}/2.$$
Now observe that the map $\textrm{B}\mathbb{Z}/3\longrightarrow
{*}$ is a $\mathbb{F}_2$-homology equivalence, and hence
$H^n(\textrm{B}\mathbf{p3m1};\mathbb{F}_2)$ is isomorphic to
$H^n(\mathbf{P}_{\textrm{B}\mathbb{Z}/3}\textrm{B}\mathbf{p3m1};\mathbb{F}_2)$.
But as $\mathbf{p3m1}$ has 2-torsion, it has nontrivial
$\mathbb{F}_2$-cohomology in arbitrarily high degrees, and then
its $\textrm{B}\mathbb{Z}/3$-nullification does, too. So, using
universal coefficients theorem, we obtain that
$\mathbf{P}_{\textrm{B}\mathbb{Z}/3}\textrm{B}\mathbf{p3m1}$ has
infinite cohomological dimension. On the other hand, as
$\underline{\textrm{B}}\mathbf{p3m1}$ is contractible and the
rational homology of $\underline{\textrm{B}}G$ is that of
$\textrm{B}G$ for any group $G$, the fact that
$\textrm{B}\mathbb{Z}/3\longrightarrow {*}$ is a
$\mathbb{Q}$-homology equivalence implies that the rational
homology of
$\mathbf{P}_{\textrm{B}\mathbb{Z}/3}\textrm{B}\mathbf{p3m1}$ is
trivial. In particular,
$\mathbf{P}_{\textrm{B}\mathbb{Z}/3}\textrm{B}\mathbf{p3m1}$ is
not homotopy equivalent to the product
$S^2\times\textrm{B}\mathbb{Z}/2$.

To conclude, note that in this context $S^2$ should be seen as two
copies of the equilateral triangle glued along their edges, with
the action of $\mathbb{Z}/2$ swapping them. Hence, the previous
fibration is not orientable (in the sense of \cite{Switzer75},
page 344). However,
$\mathbf{P}_{\textrm{B}\mathbb{Z}/3}\textrm{B}\mathbf{p3m1}$
cannot be homotopy equivalent to the projective plane, because
this nullification has infinite cohomological dimension.

We think that the ideas developed in this section can give a lot
of information about $\textrm{B}\mathbb{Z}/p$-nullification of
classifying spaces of groups of symmetries, and we plan to
undertake in subsequent work its description for all the
crystallographic (wallpaper and hyperbolic) groups, and also other
groups of symmetries as rosette or frieze groups.

\section{The homotopy fiber of the natural map $\textrm{B}G \longrightarrow
\textrm{B}_{\mathcal{F}}\emph{G}$}\label{map} \label{comacat}

We conclude this note by describing to what extent the homotopy
fiber of the map
$f:\textrm{B}G\longrightarrow\textrm{B}_{\mathcal{F}}G$ defined at
the beginning of section 3 can be built using as pieces
classifying spaces of subgroups of $G$ that belong to the family
$\mathcal{F}$. To make this decomposition, the main tools that we
are going to use are the left homotopy Kan extension of a functor
and the Gabriel-Zisman localization. Now we will recall briefly
these definitions.

In the sequel $\mathcal{C}$ and $\mathcal{D}$ will be small
topological categories. Let $F:\mathcal{C}\rightarrow\mathcal{D}$
be a functor. If \emph{d} is an object of $\mathcal{D}$, then we
define the \emph{overcategory} $F\downarrow d$ as the category
whose objects are pairs $(c,\phi)$ such that $c$ is an object of
$\mathcal{C}$ and $\phi :F(c)\rightarrow d$ is a morphism in
$\mathcal{D}$. A morphism between two pairs $(c,\phi)$ and
$(c',\phi ')$ is given by a map $\psi :c\rightarrow c'$ in
$\mathcal{C}$ such that $\phi (F(c))=\phi '\circ F(\psi)(c')$. In
the same way, the \emph{undercategory} $d\downarrow F$ is defined
as the category whose objects are pairs $(c,\phi)$ with
$c\in\mathcal{C}$ and $\phi :c\rightarrow F(d)$ a morphism in
$\mathcal{D}$. A morphism between $\psi :(c,\phi)\rightarrow
(c',\phi ')$ is a morphism $\psi ':c\rightarrow c'$ such that
$F(\psi ')\circ\phi =\phi '$. When $F\downarrow d$ (respectively
$d\downarrow F$) is contractible for every object $d$ in
$\mathcal{D}$ we say that $F$ is \emph{left cofinal} (respectively
\emph{right cofinal}).

\begin{rem}

The overcategory and undercategory are particular cases of\break
``comma categories". For a complete study of the comma categories
in the general context of category theory see (\cite{MacLane71},
II.6).

\end{rem}

From now on, and unless explicit mention against it, we will work
in the category $\mathbf{Spaces}$ of simplicial spaces (although
most of the spaces that will appear will have a simplicial complex
structure)

Let $F:\mathcal{C}\rightarrow\mathcal{D}$ be a functor. Segal
defined another functor, induced by $F$, $$L_F
:\mathbf{Fun}(\mathcal{C},\mathbf{Spaces})\longrightarrow
\mathbf{Fun}(\mathcal{D},\mathbf{Spaces}),$$ whose value on every
$X:\mathcal{C}\rightarrow\mathbf{Spaces}$ is given by
$L_F(X)(d)={\textrm{hocolim}}_{F\downarrow d}X\circ p$, where $p$
is the projection functor $p:F\downarrow d\rightarrow\mathcal{C}$.
The functor $L_F(X)$ is called the \emph{homotopy left Kan
extension} of $X$ along $F$.

The importance of that construction comes mainly from the next
result:

\begin{thm}[Homotopy pushdown theorem]

If $F:\mathcal{C}\rightarrow\mathcal{D}$ and
$X:\mathcal{C}\rightarrow\mathbf{Spaces}$ are functors, then there
is a homotopy equivalence
$$\emph{hocolim}_{\mathcal{D}}L_F(X)\simeq
\emph{hocolim}_{\mathcal{C}}X.$$ \label{HPT}
\end{thm}

\begin{proof}

The proof is done making use of the description of the homotopy
left Kan extension as the classifying space of a category. See
(\cite{Hollender92}, 5.5).  \end{proof}

Now we will recall the classical definition of localization of a
category.

\begin{thm}[Gabriel-Zisman]

Let $\mathcal{C}$ be a category. There exists another category
$\mathcal{L}(\mathcal{C})$ and a functor
$\mathcal{C}\longrightarrow\mathcal{L}(\mathcal{C})$ such that the
following conditions hold:

\begin{itemize}

\item $\mathcal{L}$ inverts the morphisms of $\mathcal{C}.$

\item If $F:\mathcal{C}\longrightarrow\mathcal{D}$ is another
functor making the morphisms of $\mathcal{C}$ invertible, there
exists one and only one functor
$F':\mathcal{L}(\mathcal{C})\longrightarrow\mathcal{D}$ such that
$F'\circ\mathcal{L}=F.$

\end{itemize}

$\mathcal{L}(\mathcal{C})$ is called the category of fractions of
$\mathcal{C}$ or simply the localization of $\mathcal{C}.$

\end{thm}

\begin{proof}

See (\cite{Gabriel67}, chapter 1).  \end{proof}

Recall that if $X$ is a simplicial complex, the simplex category
$\mathbf{\Gamma X}$ is the category whose objects are the
simplices of $X$, and whose maps are the face maps (there are no
nontrivial degeneracies). We will assume in the rest of the
section that we will work with the model of
$\textrm{B}_{\mathcal{F}}G$ constructed in proposition
\ref{modelomodific}. In the problem we are interested,
$\mathbf{\Gamma} \mathbf{B}_{\mathcal{F}}\mathbf{G}^{op}$ will
play the role of $\mathcal{C}$, and $\mathcal{D}$ will be the
localization of
$\mathbf{\Gamma}\mathbf{B}_{\mathcal{F}}\mathbf{G}^{op}$. From now
on, we will use the model of $\textrm{B}_{\mathcal{F}}G$ given in
Proposition \ref{modelomodific}, and it is not hard to see that in
this case, the $\mathbf{\Gamma}
\mathbf{B}_{\mathcal{F}}\mathbf{G}$ is exactly the category
$\mathbf{\Gamma}$ defined there.

We have developed now all the ingredients we need, and we can give
the decomposition, that is based in the concept of ``homotopy
average", proposed by Dror-Farjoun (see \cite{Dror-Farjoun95},
chapter 9). So, consider the map $\textrm{B}G\rightarrow
\textrm{B}_{\mathcal{F}}G$; if $S$ is the functor defined in
section 2, call $\overline{S}$ the composition of $\textrm{N}(S)$
with the Borel construction $\textrm{E}G\times_G{(-)}$. According
to proposition \ref{modelomodific} and (\cite{Dwyer00}, 6.5), we
have that $\textrm{hocolim}_{\mathbf{\Gamma
}^{op}}\overline{S}\simeq \textrm{B}G$.  Now, if $\mathcal{L}$ is
the localization functor previously defined, we can consider the
left homotopy Kan extension $L_{\mathcal{L}}(\overline{S})$. The
homotopy pushdown theorem \ref{HPT} implies that we have a
homotopy equivalence
$$\textrm{hocolim}_{\mathbf{\Gamma }^{op}}\overline{S}\simeq
{\textrm{hocolim}}_{\mathcal{L}(\mathbf{\Gamma}^{op})}L_{\mathcal{L}}(\overline{S}).$$
So joining all these data we obtain a string of maps
\begin{equation} \label{eq:string}
{\textrm{hocolim}}_{\mathcal{L}({\mathbf{\Gamma}}^{op})}L_{\mathcal{L}}(\overline{S})\simeq
\textrm{hocolim}_{\mathbf{\Gamma^{op}}}\overline{S}\simeq\textrm{B}G\rightarrow
\textrm{B}_{\mathcal{F}}G.
\end{equation} that, up to homotopy equivalence, is the map
$\textrm{B}G\longrightarrow\textrm{B}_{\mathcal{F}}G$ that we are
analyzing. So, we need to describe the inverse image of a simplex
$\sigma$ of $\textrm{B}_{\mathcal{F}}G$ in
${\textrm{hocolim}}_{\mathcal{L}(\mathbf{\Gamma}^{op})}L_{\mathcal{L}}(\overline{S})$.

Consider first a simplex $\sigma\in\mathbf{\Gamma}$; it can be
identified with a chain $G/H_1\rightarrow\ldots\rightarrow G/H_n$,
and then $\textrm{N}(S)(\sigma)=G/H_1$ (as a discrete $G$-set). In
a similar way, if
$$\{G/H_{i_1}\rightarrow\ldots\rightarrow G/H_{i_k}\}\longrightarrow\{ G/H_1\rightarrow\ldots\rightarrow
G/H_n\}$$ represents a morphism in $\mathbf{\Gamma}$, its image by
$\textrm{N}(S)$ is represented by the corresponding map
$G/H_1\rightarrow G/H_{i_1}$. Now, the explicit definition over a
simplex of $\textrm{N}(S)$ is immediate.

Thus, if we identify now
$\textrm{hocolim}_{\mathbf{\Gamma}^{op}}\overline{S}$ with the
Grothendieck construction of the nerve of $\overline{S}$, every
simplex of $\textrm{hocolim}_{\mathbf{\Gamma}^{op}}\overline{S}$
is represented (\cite{Dwyer00}, 6.5) by a triple $(\tau,
\sigma_1\leq\ldots\leq\sigma_n, aH_{1})$ where $\tau$ is a simplex
of $\textrm{E}G$, $\sigma_i$ is a simplex of $\mathbf{\Gamma}$,
$g\in G$ and $\sigma_n=G/H_1\rightarrow\ldots\rightarrow G/H_n$.
Observe that $\sigma_1\leq\ldots\leq\sigma_n$ represents a simplex
in $\textrm{N}(\mathbf{\Gamma})$, which is a model for
$\textrm{B}_{\mathcal{F}}G$, as we know. Now the map
$\textrm{B}G\longrightarrow\textrm{B}_{\mathcal{F}}G$ that we are
studying can be seen as the natural map
$$\textrm{hocolim}_{\mathbf{\Gamma}^{op}}\overline{S}\longrightarrow\textrm{hocolim}_{\mathbf{\Gamma}^{op}}{*}$$
induced by the natural transformation $\overline{S}\rightarrow
{*}$. Hence the image of the previous triple under this map is the
simplex $\sigma$, and according to the definition of
$\overline{S}$ we can identify the inverse image of $\sigma$ with
the pair $(\sigma ,\overline{S}(\sigma_n))$. Observe that, because
of the definition of $\overline{S}$, the latter is the same as
$(\sigma ,\overline{S}(\sigma))$, which is in fact a model for the
classifying space of $H_1$.

It remains to identify the image of the pair $(\sigma
,\overline{S}(\sigma_n))$ under the map
$$\textrm{hocolim}_{\mathbf{\Gamma}^{op}}\overline{S}\simeq
{\textrm{hocolim}}_{\mathcal{L}(\mathbf{\Gamma}^{op})}L_{\mathcal{L}}(\overline{S})$$
induced by the localization. As before, is is not hard to show
that that every simplex of
${\textrm{hocolim}}_{\mathcal{L}(\mathbf{\Gamma}^{op})}L_{\mathcal{L}}(\overline{S})$
can be written as a pair
$(\sigma_1\rightarrow\ldots\rightarrow\sigma_n,
L_{\mathcal{L}}(\overline{S})(\sigma_1)\rightarrow\ldots\rightarrow
L_{\mathcal{L}}(\overline{S})(\sigma_n))$, where
$\sigma_i\rightarrow\sigma_{i+1}$ is a morphism in
$\mathbf{\Gamma}$ or the formal inverse of a morphism in
$\mathbf{\Gamma}$. Now, from the definition of $L_{\mathcal{L}}$
it can be deduced that the image of $(\sigma ,
\overline{S}(\sigma))$ is $(\sigma
,L_{\mathcal{L}}\overline{S}(\sigma))$. Note that in the latter we
are looking at $\sigma$ in the localized category.

So we have checked that the inverse image of $\sigma$ by the
string of equivalences (\ref{eq:string}) is exactly $(\sigma
,L_{\mathcal{L}}\overline{S}(\sigma))$. Hence, we can establish
the following

\begin{thm}

If $f:\emph{B}G\rightarrow\emph{B}_{\mathcal{F}}G$ is the map
previously defined, then
$$\mathbf{Fib  } |f|\simeq |L_{\mathcal{L}}\overline{S}(\sigma)|$$ for any simplex
$\sigma$ of $\emph{B}_{\mathcal{F}}G$. Here $|\textrm{ }|$ denotes
the realization functor, and $\mathbf{Fib } |f|$ stands for the
homotopy fiber of $|f|$. \label{deco}
\end{thm}

\begin{proof}

As the category $\mathbf{\Gamma}$ is itself constructed as a
nerve, $\textrm{N}(\mathbf{\Gamma})$ has structure of simplicial
complex. So, every point of $\textrm{N}(\mathbf{\Gamma})$ belongs
to the interior of one and only one simplex of
$\textrm{N}(\mathbf{\Gamma})$, and it is enough to verify that the
fibers of the simplices of $\textrm{N}(\mathbf{\Gamma})$ are
homotopy equivalent. On the other hand, a simplex is always
contractible, so we need only check that the homotopy type of
$|L_{\mathcal{L}}\overline{S}(\sigma)|$ does not depend on the
simplex $\sigma$ of $\textrm{B}_{\mathcal{F}}G$. We know, by the
construction of the Kan extension, that
$L_{\mathcal{L}}(\overline{S})(\sigma)={\textrm{hocolim}}_{\mathcal{L}\downarrow
\sigma}(\overline{S}\circ p)$, where $p$ is the projection
 functor $p:\mathcal{L}\downarrow \sigma\rightarrow\mathcal{C}$.
So, if $\sigma$ and $\sigma'$ are two distinct simplices of
$\textrm{B}_{\mathcal{F}}G$, it is enough to see that the
overcategories $\mathcal{L}\downarrow \sigma$ and
$\mathcal{L}\downarrow \sigma'$ are equivalent. In order to check
this, let $g:\sigma\rightarrow\sigma'$ be a morphism in
$\mathcal{L}(\mathbf{\Gamma}^{op})$, that always exist because
$\textrm{B}_{\mathcal{F}}G$ is connected. In these conditions, $g$
induces a natural transformation
$$\mathbf{T}_g:\mathcal{L}\downarrow \sigma\longrightarrow
\mathcal{L}\downarrow \sigma'$$ that sends every object $(\tau
,\phi)$ of $\mathcal{L}\downarrow \sigma$ to $(\tau
,g\circ\phi)\in\mathcal{L}\downarrow \sigma'$ and the morphisms to
the obvious ones. But the morphism $g$ is invertible (because it
is a morphism in the localized category), and clearly the natural
transformations $\mathbf{T}_g$ and $\mathbf{T}_{g^{-1}}$ are
inverses one of each other. In other words, the two overcategories
are equivalent, and the corresponding homotopy colimits have the
same homotopy type. So we are done.  \end{proof}

The following corollary is immediate:

\begin{cor}

The homotopy fiber $\mathbf{Fib}|f|$ has the homotopy type of\linebreak
$|{\emph{hocolim}}_{\mathcal{L}\downarrow \sigma}$
$(\overline{S}\circ p)|$, and in particular it is a homotopy
colimit of classifying spaces of groups of $\mathcal{F}$ over a
contractible category.

\end{cor}

We conclude this section by proving that the nerves of the two
overcategories that appear in the proof of the previous theorem
are contractible. We think that this question can have independent
interest, and we would like to point out that, although the result
seems to be known (see \cite{Dror-Farjoun95}, 9.E.3), we have
found no proof in the literature, so we give this one.

\begin{prop}

Let $X$ be a simplicial complex, and let
$\mathcal{L}:\mathbf{\Gamma
X}\longrightarrow\mathcal{L}(\mathbf{\Gamma X})$ be the
Gabriel-Zisman localization functor, where $\mathbf{\Gamma X}$ is
the simplex category of $X$. Then for every simplex $\sigma\in X$
the overcategory $\mathcal{L}\downarrow\sigma$ is contractible.

\end{prop}

\begin{proof}

The idea of the proof is to build, for every simplex $\sigma\in X$
a homotopy between the identity map
$\textrm{Id}_{|\mathbf{N}(\mathcal{L}\downarrow\sigma)|}$ and a
constant map. In order to do this, we will prove firstly the
existence of a sequence of endofunctors
$$\{F_n\}:\mathcal{L}\downarrow\sigma\longrightarrow\mathcal{L}\downarrow\sigma
$$ for every $n\geq 0$ such that $F_0=\textrm{Id}$ and for every
$(\tau,\phi)\in\mathcal{L}\downarrow\sigma$ there exists a natural
number $n_{(\tau,\phi)}$ in such a way that
$F_m((\tau,\phi))=(\sigma,\textrm{Id})$ for every $m\geq
n_{(\tau,\phi)}$.

In the sequel the maps in $\mathbf{\Gamma X}$ and their images in
$\mathcal{L}(\mathbf{\Gamma X})$ will be denoted indistinctly by
$i_{\alpha}$, where $\alpha$ will be an appropriate subindex. The
inverse of $i_{\alpha}$ in the localized category will be called
$j_{\alpha}$.

It is plain from the definition of the localization functor that
every element $(\tau,\phi)$ of $\mathcal{L}\downarrow\sigma$
admits a unique expression of the form $(\tau, j_n\circ
i_{n-1}\circ\ldots\circ j_2\circ i_1)$, where we allow that $j_n$
or $i_1$ can be the identity (but no one of the other maps that
appear), $j_{t-1}\neq i^{-1}_t\neq j_{t+1}$ for every $t$.

So, we begin with $F_0=\textrm{Id}$. Let us define the functor
$$F_1:\mathcal{L}\downarrow\sigma\longrightarrow\mathcal{L}\downarrow\sigma .$$
If $(\tau, j_n\circ i_{n-1}\circ\ldots\circ j_2\circ i_1)$ is an
element of the overcategory, then we say $F_1((\tau, j_n\circ
i_{n-1}\circ\ldots\circ j_2\circ i_1))=(i_1(\tau), j_n\circ
i_{n-1}\circ\ldots\circ j_2)$, and the map induced by a face map
will be sent to the identity map between the images. It is easy to
see that this functor is well-defined.

Now,
$F_2:\mathcal{L}\downarrow\sigma\longrightarrow\mathcal{L}\downarrow\sigma$
will be defined as $F_2((\tau, j_n\circ i_{n-1}\circ\ldots\circ
j_2\circ i_1))=(j_2^{-1}\circ i_1(\tau), j_n\circ
i_{n-1}\circ\ldots\circ i_3)$. Observe that this is well-defined
because the localization functor is, in this case, bijective over
the objects. Again, the image of every morphism by $F_2$ will be
the identity. It is clear again that this is a functor.

In an analogous way, we can define, for $m$ odd, $F_m((\tau,
j_n\circ i_{n-1}\circ\ldots\circ j_2\circ i_1))=(i_{m}\circ
j_{m-1}\circ i_1(\tau), j_n\circ i_{n-1}\circ\ldots j_{m+1})$, and
for $m$ even, $F_m((\tau, j_n\circ i_{n-1}\circ\ldots\circ
j_2\circ i_1))=(j_{m}\circ i_{m-1}\circ i_1(\tau), j_n\circ
i_{n-1}\circ\ldots i_{m+1})$, and the image sends every morphism
to the identity map. This is the sequence we were looking for.

Our next goal will be to relate all these functors by natural
transformations, in order to obtain the desired homotopy.

Let $m\geq 0$ be again a natural number. First we will define the
transformation $T_{2m}:F_{2m}\longrightarrow F_{2m+1}$. If $(\tau,
j_n\circ i_{n-1}\circ\ldots\circ j_2\circ i_1)$ is an object of
the overcategory, we define the map $F_{2m}((\tau, j_n\circ
i_{n-1}\circ\ldots\circ j_2\circ i_1))\longrightarrow
F_{2m+1}((\tau, j_n\circ i_{n-1}\circ\ldots\circ j_2\circ i_1))$
as the obvious map induced by $$i_{2m}:j_{2m-1}\circ\ldots\circ
i_1(\tau)\longrightarrow i_{2m}\circ j_{2m-1}\circ\ldots\circ
i_1(\tau).$$
On the other hand, we define, for every $m\geq 1$, the natural
transformation $T_{2m-1}:F_{2m}\longrightarrow F_{2m-1}$ in the
following way: $F_{2m}((\tau, j_n\circ i_{n-1}\circ\ldots\circ
j_2\circ i_1))\longrightarrow F_{2m-1}((\tau, j_n\circ
i_{n-1}\circ\ldots\circ j_2\circ i_1))$ is the map induced by
$$i_{2m-1}:j_{2m-1}\circ\ldots\circ i_1(\tau)\longrightarrow
i_{2m-2}\circ j_{2m-3}\circ\ldots\circ i_1(\tau).$$ Recall the
fact that, by definition of the $j$'s, $j^{-1}$ represents a
morphism in $\mathbf{\Gamma X}.$

By the previous arguments we have defined a string of natural
transformations
$$\textrm{Id}=F_0\stackrel{T_0}{\longrightarrow}F_1\stackrel{T_1}{\longleftarrow}
F_2\stackrel{T_2}{\longrightarrow}F_3\stackrel{T_3}{\longleftarrow}\ldots$$
Before we continue, we shall do a couple of remarks.

\begin{itemize}

\item It is known (\cite{Dwyer00}, I.5) the functors $F_n$ define
simplicial maps from nerves
$$\textrm{N}(F_n):\textrm{N}(\mathcal{L}\downarrow\sigma)\longrightarrow
\textrm{N}(\mathcal{L}\downarrow\sigma)$$ which, the same way,
define maps $|f_n|$ from the realization of the nerve to itself.
The fact that $F_n$ is always related to $F_{n+1}$ by a natural
transformations means that $f_n$ is simplicially homotopic to
$f_{n+1}$, and, in addition, $|f_n|$ is homotopic to $|f_{n+1}|$.
The crucial point here is the homotopies between the realization
of the maps are first defined over the vertices of the nerve of
$\mathcal{L}\downarrow\sigma\times I$ (with the usual simplicial
structure of the product) and then extended by linearity to all
the complex. We will use this fact later.

\item Let $(\tau, j_n\circ\ldots\circ i_{1})$ be an object of the
overcategory. From the definitions of the functors $F_i$ we can
deduce that $F_n\circ\ldots\circ F_1 ((\tau, j_n\circ\ldots\circ
i_{1}))=(\sigma,\textrm{Id})$. So, as the chain of maps
$j_n\circ\ldots\circ i_{1}$ is always finite, we can say that for
every $(\tau,\phi)\in (\mathcal{L}\downarrow\sigma)$ there exists
a minimal natural number $n_{(\tau,\phi)}$ such that
$F_{n_{\tau,\phi}}\circ\ldots\circ F_1 ((\tau,
\phi))=(\sigma,\textrm{Id})$. At the level of nerves, we are
saying that for every vertex
$v\in\textrm{N}(\mathcal{L}\downarrow\sigma)$ there exists $n_v$
such that $f_{n_v}\circ\ldots\circ
f_1(v)=\textrm{N}(\sigma,\textrm{Id})$.

\end{itemize}

For $n$ even, let us call $H_n$ the simplicial homotopy induced by
the transformation $T_n$. If $n$ is odd, we call $H'_{n-1}$ the
homotopy induced by $T_n$ between $f_n$ and $f_{n-1}$, and put
$H_{n-1}(x,t)=H'_{n-1}(x,1-t)$, the homotopy that begins in
$f_{n-1}$ and ends at $f_n$.

Now we are prepared to define the homotopy between the identity
and the constant map from the realization to itself with value
$|\textrm{N}(\sigma,\textrm{Id})|$ (in the rest we will call this
element $*$). So, consider a vertex $v\in
\textrm{N}(\mathcal{L}\downarrow\sigma)$. We define a map
$H:|\textrm{N}(\mathcal{L}\downarrow)\sigma|\times
I\longrightarrow |\textrm{N}(\mathcal{L}\downarrow)\sigma|$ by
$$ H(v,t)=\left\{ \begin{array}{ll} |H_0|(v,n_vt) & \textrm{if
}t\in [0,\frac{1}{n_v}]\\ |H_1|(v,n_vt-1) & \textrm{if }t\in
[\frac{1}{n_v},\frac{2}{n_v}]\\ \vdots & \vdots \\
|H_{n-1}|(v,n_vt-(n-1)) & \textrm{if }t\in [\frac{n_v-1}{n_v},1]
\end{array} \right.
$$
The map $H$ defined in this way lineally extends to all of
$|\textrm{N}(\mathcal{L}\downarrow\sigma)|$. Let us see that $H$
is the desired map.

\begin{enumerate}

\item If $v$ is a vertex of $\textrm{N}(\mathcal{L}\downarrow\sigma)$,
$H(v,0)=H_0(v,0)=v$. In the same way, $H(v,1)=H_{n_v}(v,1)=*$. As
$|H|_i$ is defined by linear extension for every $i$ and the same
happens with $H$, the previous equalities hold for every point of
the complex.

\item $H$ is continuous with respect to $t$ because the homotopies $|H_i|$
are, and $|H_j(x,1)|=f_{j+1}(x)=|H_{j+1}(x,0)|$ for every
$x\in|\textrm{N}(\mathcal{L}\downarrow\sigma)|$.

\item Finally, $H$ is continuous respect the first component
because it is defined by linear extension of a map defined on the
vertices of a simplicial complex.

\end{enumerate}

These three statements prove that $H$ is the homotopy between the
identity and the constant map we were looking for. So,
$\mathcal{L}\downarrow\sigma$ is contractible.  \end{proof}

\Addresses\recd

\end{document}

%% file: agtout.tex
\def\ifplaintex{\expandafter\ifx\csname documentclass\endcsname\relax}

\def\gtp{{\mathsurround=0pt\it $\cal G\mskip-2mu$eometry \&\ 
$\cal T\!\!$opology $\cal P\!$ublications}}  % GT publications

\def\recd{{\small Received:\qua\receiveddate\ifx\reviseddate\relax
\else\qquad Revised:\qua\reviseddate\fi\par}} 

\def\lognumber#1{\def\thelognumber{#1}}
\def\volumenumber#1{\def\thevolumenumber{#1}}
\def\volumeyear#1{\def\thevolumeyear{#1}}
\def\papernumber#1{\def\thepapernumber{#1}}
\def\pagenumbers#1#2{\def\startpage{#1}\def\finishpage{#2}}
\def\published#1{\def\publishdate{#1}}

\def\received#1{\def\receiveddate{#1}}
\def\revised#1{\def\reviseddate{#1}}
\def\accepted#1{\def\accepteddate{#1}}

\def\asciiauthors#1{\def\theasciiauthors{#1}}
\def\asciiaddress#1{\def\theasciiaddress{#1}}

\def\asciiurl#1{\def\theasciiurl{#1}}
\def\coverauthors#1{\def\thecoverauthors{#1}}
\long\def\asciiabstract#1{\long\def\theasciiabstract{#1}}

\let\\\par\let\thelognumber\relax\let\thevolumenumber\relax
\let\thepapernumber\relax\let\thevolumeyear\relax\let\startpage\relax
\let\finishpage\relax\let\publishdate\relax\let\receiveddate\relax
\let\reviseddate\relax\let\accepteddate\relax\let\theasciititle\relax
\let\theasciiauthors\relax\let\theasciiaddress\relax
\let\theasciiabstract\relax

\let\thecoverauthors\relax\let\theasciiemail\relax
\let\theasciiurl\relax

\ifplaintex
\font\logobig=cmssbx10 scaled 3836
\font\logomed=cmssbx10 scaled 2557
\else
\font\logobig=cmssbx10 scaled 4200
\font\logomed=cmssbx10 scaled 2800
\fi

\long\def\makeagttitle{   %%% start of definition of \makeagttitle
\count0=\startpage
\agt\hfill      %   Journal title (top left) 
\hbox to 45truept{\vbox to 0pt{\vglue -13truept{\logomed A\kern -.37em{\logobig 
T}\kern -.38em G}\vss}\hss}
\break
{\small Volume \thevolumenumber\ (\thevolumeyear)
\startpage--\finishpage\nl
Published: \publishdate}

\vglue .25truein

{\parskip=0pt\leftskip 0pt plus
1fil\def\\{\par\smallskip}{\Large\bf\thetitle}\par\medskip} \vglue
0.05truein

{\parskip=0pt\leftskip 0pt plus 1fil\def\\{\par}{\sc\theauthors}
\par\medskip}%
 
\vglue 0.03truein 

{\small\leftskip 25truept\rightskip 25truept{\bf Abstract}\stdspace\theabstract

{\bf AMS Classification}\stdspace\theprimaryclass
\ifx\thesecondaryclass\relax\else; \thesecondaryclass\fi\par
{\bf Keywords}\stdspace \thekeywords\par}\vglue 7truept

}   %%%% end of definition of \makeagttitle

\ifplaintex
\font\phead=cmsl9 scaled 950
\font\pnum=cmbx10 scaled 913
\font\pfoot=cmsl9 scaled 950
\headline{\vbox to 0pt{\vskip -4.5mm\line{\small\phead\ifnum
\count0=\startpage ISSN 1472-2739 (on-line) 1472-2747 (printed)
\hfill {\pnum\folio}\else\ifodd\count0\def\\{ }% 
\ifx\theshorttitle\relax\thetitle\else\theshorttitle\fi\hfill{\pnum\folio}
\else\def\\{ and }{\pnum\folio}\hfill\ifx\theshortauthors\relax\theauthors
\else\theshortauthors\fi\fi\fi}\vss}}
\footline{\vbox to 0pt{\vglue 0mm\line{\small\pfoot\ifnum\count0=\startpage
\copyright\ \gtp\hfill\else
\agt, Volume \thevolumenumber\ (\thevolumeyear)\hfill\fi}\vss}}
\else
\font=cmsl9 scaled 1050
\font\lnum=cmbx10 
\font=cmsl9 scaled 1050
\makeatletter
\def\@oddhead{{\small\lhead\ifnum\count0=\startpage ISSN 1472-2739 
(on-line) 1472-2747 (printed)\hfill {\lnum\number\count0}\else\ifodd\count0
\def\\{ }\ifx\theshorttitle\relax \thetitle \else\theshorttitle\fi\hfill
{\lnum\number\count0}\else\def\\{ and }{\lnum\number\count0}
\hfill\ifx\theshortauthors\relax 
\theauthors\else\theshortauthors\fi\fi\fi}}\def\@evenhead{\@oddhead}
\def\@oddfoot{\small\lfoot\ifnum\count0=\startpage\copyright\ \gtp\hfill\else
\agt, Volume \thevolumenumber\ (\thevolumeyear)\hfill\fi}
\def\@evenfoot{\@oddfoot}
\makeatother
\fi
\let\maketitlepage\makeagttitle
\let\makeshorttitle\maketitlepage
\let\maketitle\maketitlepage

\newwrite\gtoutfile
\long\gdef\makeheadfile{  %%% start of definition of \makeheadfile
{\def\\{, }\def\s{ }
\immediate\openout\gtoutfile head.xxx
\immediate\write\gtoutfile{Proxy-for: \ifx\theasciiauthors\relax
\theauthors\else\theasciiauthors\fi\s<\ifx\theasciiemail\relax\theemail\else\theasciiemail\fi>}
\immediate\write\gtoutfile{\noexpand\\}
\immediate\write\gtoutfile{Authors: \ifx\theasciiauthors\relax
\theauthors\else\theasciiauthors\fi}
{\def\\{ }\immediate\write\gtoutfile{Title: \ifx\theasciititle\relax
\thetitle\else\theasciititle\fi}}
\immediate\write\gtoutfile{Subj-class: GT or SG, GR etc}
\immediate\write\gtoutfile{MSC-class: \theprimaryclass\ifx\thesecondaryclass\relax\else, \thesecondaryclass\fi}
\immediate\write\gtoutfile{Journal-ref: Algebr. Geom. Topol. \thevolumenumber\s
(\thevolumeyear) \startpage-\finishpage}
\immediate\write\gtoutfile{Comments: Published by Algebraic and
Geometric Topology at}
\immediate\write\gtoutfile{\s\s\s  http://www.maths.warwick.ac.uk/agt/AGTVol\thevolumenumber/agt-\thevolumenumber-\thepapernumber.abs.html}
\immediate\write\gtoutfile{\noexpand\\}
\immediate\write\gtoutfile{}
\ifx\theasciiabstract\relax
\immediate\write\gtoutfile{\theabstract}\else
\immediate\write\gtoutfile{\theasciiabstract}\fi
\immediate\write\gtoutfile{}
\immediate\write\gtoutfile{\noexpand\\}
\immediate\write\gtoutfile{}
\immediate\closeout\gtoutfile}}  %%% end of definition of \makeheadfile

\def\maketitlepage{\makeagttitle\makeheadfile}
\let\makeshorttitle\maketitlepage
\let\maketitle\maketitlepage